\documentclass[final,leqno]{siamltex}

\usepackage{amsmath}
\usepackage{amssymb}
\usepackage[usenames, dvipsnames]{color}
\usepackage[ruled, vlined]{algorithm2e}
\usepackage{bigints}
\usepackage{esint}
\usepackage{bm}
\usepackage{color}
\usepackage{comment}
\usepackage{graphicx}
\usepackage{subfig, overpic}
\usepackage[capitalize,nameinlink]{cleveref} \crefname{equation}{}{Equations}
\Crefname{equation}{Eq.}{Equations}

\newcommand{\sampleSetSize}{p}
\newcommand{\mixingCoeff}{m}
\newcommand{\epsilonMul}{\eta}
\newcommand{\problemDim}{n}
\newcommand{\reducedBasisDim}{n_{rb}}
\newcommand{\skel}[1]{\widehat{#1}}
\newcommand{\paramSet}{\Omega}
\newcommand{\param}{\omega}
\newcommand{\domain}{\mathcal{D}}
\newcommand{\operatorNP}{\mat{L}}
\newcommand{\operator}[1]{\operatorNP(#1)}
\newcommand{\operatorCoarseNP}{\operatorNP_C}
\newcommand{\operatorCoarse}[1]{\operatorCoarseNP(#1)}
\newcommand{\mat}[1]{\mathsf{#1}}
\newcommand{\ve}[1]{\mathsf{#1}}
\newcommand{\solution}[1]{\ve{u}(#1)}
\newcommand{\solutions}{\mat{S}}

\newcommand{\solutionsCoarse}{\solutions_C}
\newcommand{\solutionCoarse}[1]{\ve{u}_C(#1)}
\newcommand{\source}{\ve{f}(\param)}
\newcommand{\sourceRB}{\ve{f}_{rb}}
\newcommand{\reducedBasis}{\mat{Q}}
\newcommand{\permutation}{\pi}
\newcommand{\mixing}{\mat{M}}
\newcommand{\operatorAffineOffset}{\mat{A}}
\newcommand{\operatorAffineParamNP}{\mat{B}}
\newcommand{\sourceCoarse}[1]{\ve{f}_C(#1)}
\newcommand{\operatorRBNP}{\operatorNP_{rb}}
\newcommand{\operatorAffineParam}[1]{\operatorAffineParamNP(#1)}
\newcommand{\skelCount}{s}
\newcommand\second{\,\text{sec}}

\newcommand{\rbApprox}{\ve{u}_{rb}}
\newcommand{\rbApproxParam}[1]{\rbApprox(#1)}
\newcommand{\rbCoeff}{\ve{v}}
\newcommand{\rbCoeffParam}[1]{\rbCoeff(#1)}
\newcommand\T{\mathrm{T}}
\newcommand\dd{\,\mathrm{d}} 
\SetKwInput{KwData}{Input}
\SetKwInput{KwResult}{Output}

\title{Coarse-Proxy Reduced Basis Methods for Integral Equations} 

\author{
  Philip A. Etter  \thanks{ICME, Stanford University, Stanford, CA 94305 ({\tt paetter@stanford.edu}).}
  \and 
  Yuwei Fan  \thanks{Department of Mathematics, Stanford University, Stanford, CA 94305 ({\tt ywfan@stanford.edu}).}
  \and
  Lexing Ying  \thanks{Department of Mathematics and ICME, Stanford University, Stanford, CA 94305 ({\tt lexing@stanford.edu}).}
}

\begin{document}
\maketitle

\begin{abstract}
In this paper, we introduce a new reduced basis methodology for accelerating the computation of
large parameterized systems of high-fidelity integral equations. Core to our methodology is the use
of \emph{coarse-proxy} models (i.e., lower resolution variants of the underlying high-fidelity
equations) to identify important samples in the parameter space from which a high quality reduced
basis is then constructed. Unlike the more traditional POD or greedy methods for reduced basis
construction, our methodology has the benefit of being both easy to implement and embarrassingly
parallel. We apply our methodology to the under-served area of integral equations, where the density
of the underlying integral operators has traditionally made reduced basis methods difficult to
apply. To handle this difficulty, we introduce an operator interpolation technique, based on random
sub-sampling, that is aimed specifically at integral operators. To demonstrate the effectiveness of our
techniques, we present two numerical case studies, based on the Radiative Transport Equation and a
boundary integral formation of the Laplace Equation  respectively, where our methodology provides a
significant improvement in performance over the underlying high-fidelity models for a wide
range of error tolerances. Moreover, we demonstrate that for these problems, as the coarse-proxy
selection threshold is made more aggressive, the approximation error of our method decreases at an
approximately linear rate.
\end{abstract}

\begin{keywords}
  Coarse-proxy; Reduced Basis Method; Model Order Reduction; Skeleton Extraction; Integral Equations 
\end{keywords}

\begin{AMS}
    65C30, 45A05, 65R20
\end{AMS}

\pagestyle{myheadings}
\thispagestyle{plain}

\section{Introduction}

Across virtually all areas of science and engineering, physical simulation has become an absolutely indispensable tool for the advancement of knowledge and the design of industrial products. However, as with any tool, there are always practical caveats. In particular, high-fidelity simulations often require tremendous computational resources and time to execute. This computational cost often precludes high-fidelity simulations from being used in many important problems, such as uncertainty quantification or Bayesian inference, that require not just one, but many queries to the underlying computational model. Making these many-query problems tractable often requires fast approximation techniques to mitigate the sheer computational cost of multiple queries to the underlying (full-order) model.

One such class of approximation techniques is reduced order models (ROMs). Reduced order models
typically operate in two stages. First, there is a computationally expensive \textit{offline} stage (i.e., \textit{training} stage),
wherein the ROM is trained on a collection of solutions to the full-order model (FOM). In many
cases, this entails finding a basis for a low-dimensional linear subspace which captures solutions
to the full-order model (i.e., a reduced basis). Once this offline stage is complete, the
reduced-order model can be deployed in an \textit{online} stage (i.e., \textit{test} stage), where these methods can compute
fast approximations to new problem instances by exploiting the problem structure learned during the
offline phase. For example, one can project the new problem instance onto a set of reduced basis and
solve a low-dimensional reduced problem instead of the high-dimensional full-order problem. We refer
the reader to \cite{hesthaven2016certified, benner2013survey} and references therein for a more
thorough overview of this topic. 

In this paper, we deal specifically with the class of ROM techniques falling under the reduced basis
method (RBM) \cite{rozza2007reduced, patera2007reduced}. The groundwork for the reduced basis method
was set in the late 1970s with work on the approximation for nonlinear structure analysis
\cite{almroth1978automatic, nagy1979modal, noor1980reduced}, particularly for beams and arches. This
groundwork later evolved into a more general framework for parameterized differential equations
\cite{gunzburger2012finite, peterson1989reduced}, with a corresponding swath of mathematical
analyses of the approximation error of the method \cite{rheinboldt1993theory, barrett1995reduced,
fink1983error, porsching1985estimation}. These nascent methods typically involved finding a low
dimensional approximation space around a parameter of interest --- thereby making them local
approximation methods. Later, this line of inquiry evolved into finding a global approximation space
constructed from a sparse set of sampled solutions to the full-order model
\cite{balmes1996parametric, ito1998reduced}. More recently, the first theoretical a priori
convergence guarantee was proved and numerically confirmed in \cite{maday2002priori}. This
demonstrated the potential of reduced basis methods as a robust approximation for parameterized
partial differential equations.

However, while these techniques are well-established for ordinary and partial differential
equations, there has been relatively little work done in the regime of model order reduction for
integral equations. The current most notable contributions in this underserved area are taiylored
specifically to boundary element formulations of the electric field equations
\cite{fares2011reduced, hesthaven2012certified, phillips1996efficient, ganesh2012reduced}. The chief
factor that contributes to this research gap is likely the difficulties that come from the operators
that arise from discretizing integral equations, which are typically dense. This operator density
precludes one from assembling the operators outright, which limits the applicability of many
existing model order reduction techniques, in part because even sampling a single entry of the
problem residual takes time on the order of the problem size. Regardless, this gap in the literature
is unfortunate, as integral equations have many desirable properties over their differential
counterparts. Integral equations are often better conditioned than differential equations, and many
important physical models, such as electromagnetism and radiative transport, are amenable to special
integral formulations with desirable properties (e.g., boundary integral formulation).

\section{Problem Statement}
The goal of this paper is to solve parameterized integral equations of the form
\begin{equation} \label{basic_eq}
  \operator{\param} \solution{\param} = \source, \qquad \param \in \paramSet_\infty,
\end{equation}
where $\operator{\param} \in \mathbb{R}^{\problemDim \times \problemDim}$ denotes a (dense) linear
elliptic integral operator, $\source \in \mathbb{R}^{\problemDim}$ denotes a source term, and
$\param$ are parameters taken from some sample space $\paramSet_\infty$. The underlying sample space
$\paramSet_\infty$ is typically continuous with respect to $\param$, so in this paper we concern
ourselves with a discrete subset $\paramSet$ of $\paramSet_\infty$, appropriately spaced so that
every point in $\paramSet_\infty$ is relatively close to a proxy or set of proxies in $\paramSet$.
Approximate solutions to equations whose parameters come from outside of $\paramSet$ can then be
formed via interpolation. Throughout this paper, we represent $\paramSet$ as a set
\begin{equation}
  \paramSet \equiv \left\{\param_1, \param_2, \dots, \param_\sampleSetSize \right\},
\end{equation} 
whose elements $\param_i$ denote the samples for which we would like to solve the integral equation
\cref{basic_eq}. 

If the parameter $\param$ wildly changes the underlying problem, then it is difficult to perform
this task more efficiently then simply solving all of the equations \cref{basic_eq}. However, in
many real-world contexts, the dependence on the parameter $\param$ is such that the solutions
$\solution{\param}$ form a space that is approximately low dimensional. In this case, the solutions
$\solution{\param}$ can be well represented by a few appropriately chosen degrees of freedom. The
goal of reduced basis methods is to extract these relevant degrees of freedom and use them to
accelerate the computation of the solutions $\solution{\param}$. 

Therefore, we ultimately want to find a small orthogonal basis matrix $\reducedBasis \in
\mathbb{R}^{\problemDim \times \reducedBasisDim}$, where $\reducedBasisDim \ll \problemDim$, whose
columns approximately capture the solution set 
\begin{equation}
    \solutions \equiv \begin{bmatrix} 
        \solution{\param_1} & \solution{\param_2} & \dots &
        \solution{\param_\sampleSetSize} 
    \end{bmatrix}.
\end{equation} 
Once the basis matrix $\reducedBasis$ is given, the solution $\solution{\param}$ can be approximated
by $\reducedBasis\reducedBasis^{\T}\solution{\param}$, and then applying the Galerkin projection on the 
system \eqref{basic_eq} yields
\begin{equation} \label{rb_system}
  [\reducedBasis^\T \operator{\param} \reducedBasis]  [\reducedBasis^\T \solution{\param}]
  \approx [\reducedBasis^\T \source].
\end{equation}
Since the dimension $\reducedBasisDim$ is much less than the dimension $\problemDim$ of the original
system \cref{basic_eq}, the projected system \cref{rb_system} provides us with an inexpensive way
of computing approximations to the solutions $\solution{\param}$. First, one solves for the quantity
$\reducedBasis^\T \solution{\param}$ in the $\reducedBasisDim \times \reducedBasisDim$ projected
system \cref{rb_system}. Afterwards, applying the matrix $\reducedBasis$ to the result
$\reducedBasis^\T \solution{\param}$ gives an approximation of the true
solution $\solution{\param}$.

\subsection{Main difficulties}
In the procedure of solving \eqref{rb_system}, there are two practical difficulties which arise:
\begin{enumerate}
  \item \emph{Assembling the reduced basis $\reducedBasis$ efficiently. (Offline).} There are a
    number of existing methods for constructing the basis $\reducedBasis$. Unfortunately, they are
    typically either computationally expensive or difficult to implement. One can perform a proper
    orthogonal decomposition (POD) of solutions to the full-order model to obtain such a basis
    $\reducedBasis$ \cite{hesthaven2016certified}. However, this requires a significant number of
    solves to the underlying full-order model. There are also greedy methods
    \cite{boyaval2010reduced, maday2002priori, maday2002global, rozza2007reduced,
    veroy2003posteriori, binev2011convergence}, which sequentially build up a reduced basis by
    repeatedly selecting the solution $\solution{\param_i}$ which would yield the greatest reduction in
    error, according to some error estimator. Unfortunately, the implementation and construction of
    error estimators are very involved. Moreover, the inherently sequential nature of greedy
    selection algorithms means that they are very difficult to parallelize. 
  \item \emph{Assembling the projected operator $\reducedBasis^\T \operator{\param} \reducedBasis$
    efficiently. (Online).} Since $\operator{\param}$ is an operator, it is usually too
    computationally expensive to assemble the whole operator $\operator{\param}$ explicitly. 
                This means that a reduced basis method must provide a cheap way of constructing the reduced
    operators $\reducedBasis^\T \operator{\param} \reducedBasis$ without ever explicitly assembling
    their full-order counterparts $\operator{\param}$. 
\end{enumerate}

\subsection{Contribution}
In this paper, we present a novel reduced basis approach to integral equations that has the benefit of being both general-purpose and easy to implement. Our contributions are twofold: we first present a novel method for \textit{efficient selection of training samples}. We use this selection scheme to address the first issue of assembling the reduced basis $\reducedBasis$ efficiently. Next, we present a \textit{simple interpolation technique for assembling reduced operators}. We use this technique to address the second difficultly of assembling the projected operator $\reducedBasis^\T \operator{\param} \reducedBasis$ efficiently. 

\begin{figure}[ht]
  \centering
  \includegraphics[width=0.9\textwidth]{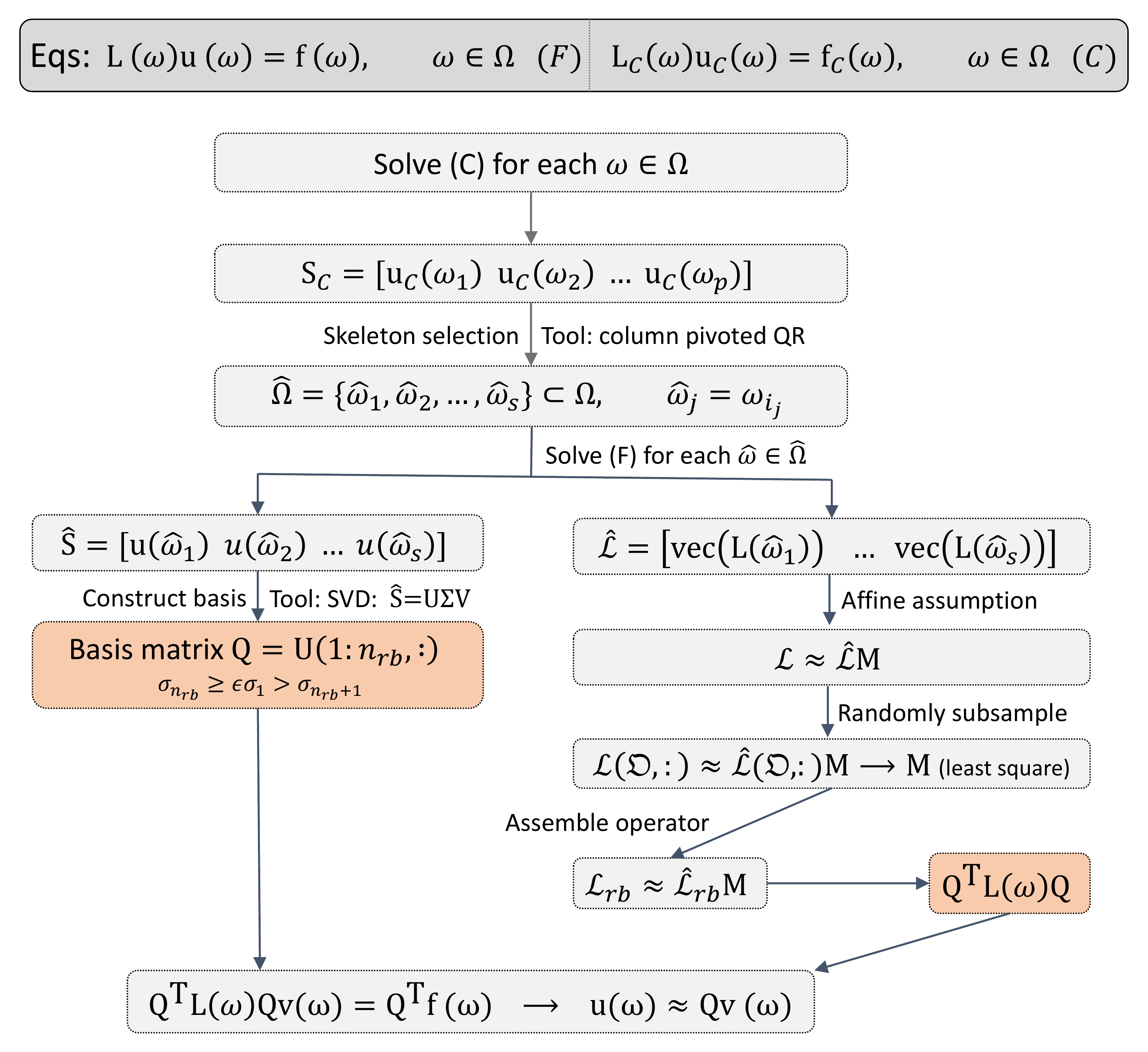}
  \caption{\label{fig:framework} Diagram of the reduced method for the integral equation.}
\end{figure}

The combination of these two techniques forms the core of our \textit{coarse-proxy} reduced basis method, diagrammed in full in \cref{fig:framework}. Our method provides a model order reduction framework for general linear integral equation problems that addresses both of the above issues without the aforementioned pitfalls of
existing methods. In our numerical experiments, we apply our method to two examples, the radiative transport equation and the boundary integral formulation of the Laplace equation---and demonstrate that on both of these problems, our methods allow for significant improvements in performance over a naive solve of all elements of the parameter space. 

For the aforementioned \textit{efficient selection of training samples}, we propose a novel method of constructing the reduced basis
$\reducedBasis$ by leveraging a coarse-proxy model to identify a set of important parameters $\skel{\param}_1, ...,
\skel{\param}_s$ in the sample space $\paramSet$, where $s\ll p$. As an example, one can use an
inexpensive low-resolution model to identify which parameters $\param$ would be important in the construction of a
reduced basis, and we only solve the full-order equations \eqref{basic_eq} for these important parameters. The method has the
desirable property of being embarrassingly parallelizable.

For the aforementioned \textit{interpolation technique for assembling reduced operators}, we propose assembling the operators $\reducedBasis^\T \operator{\param}
\reducedBasis$ by levering the power of random sampling. To be more precise, we draw random samples of the
operators $\operator{\param}$ and then use these samples to linearly interpolate between a subset of basis operators
$\reducedBasis^\T \operator{\skel{\param}_1} \reducedBasis$, ..., $\reducedBasis^\T
\operator{\skel{\param}_s} \reducedBasis$ to approximately reconstruct $\reducedBasis^\T \operator{\param}
\reducedBasis$.  This method is similar to the matrix gappy POD technique
proposed in \cite{carlberg2015preserving} by Carlberg et al, but slightly different because we cannot
afford to orthogonalize operators.

The details of the proposed method are discussed in \cref{sec:framework}, and numerical tests are
presented in \cref{sec:num}.

\section{Framework Details}\label{sec:framework}

Our framework for solving problems of the form in \eqref{basic_eq} is based on the idea of using
an inexpensive coarse-proxy model to extract the important solutions from the solution set $\solutions$. This
model can be, for example, the original fine problem, but at a much lower resolution. 
Alternately, one can use a sparse basis of wavelets. 

We write this coarse-proxy model as 
\begin{equation} 
  \operatorCoarse{\param} \solutionCoarse{\param} = \sourceCoarse{\param}.  
\end{equation}
where $\operatorCoarse{\param} \in \mathbb{R}^{n_C \times n_C}$ is the coarse analogue of the
operator $\operator{\param}$, and $\solutionCoarse{\param} \in \mathbb{R}^{n_C}$ and
$\sourceCoarse{\param} \in \mathbb{R}^{n_C}$ are the coarse-proxy solution and coarse-proxy source term
respectively. One should choose this coarse-proxy model so that it is inexpensive to evaluate (i.e., $n_C^2 \ll
\problemDim^2$). But, as long as the solutions $\solutionCoarse{\param}$ of the coarse-proxy model can
approximately capture the important features of their fine counterparts $\solution{\param}$, the
particular choice of coarse-proxy model is not especially relevant. However, one must
still exercise the appropriate caution. For example, if the solutions $\solution{\param}$ contain important high
frequency content, one should not expect that solving the problem on a coarse grid will provide a
good coarse-proxy model. 

\paragraph{Notation}
We use MATLAB notation to denote submatrices, i.e., if $\mat{M} \in \mathbb{R}^{n \times m}$, then for $A \subseteq \{1, ..., n\}$ and $B \subseteq \{1, ..., m\}$, $\mat{M}(A, B) \in \mathbb{R}^{|A| \times |B|}$ denotes the submatrix of $\mat{M}$ formed with rows $A$ and columns $B$. In the case where either $A = \{1, ..., n\}$ or $B = \{1, ..., m\}$, we use the shorthand $\mat{M}(:, B) \in \mathbb{R}^{n \times |B|}$ or $\mat{M}(A, :) \in \mathbb{R}^{|A| \times m}$, respectively. The same notation also applies to vectors.

\subsection{Skeleton Extraction} \label{basis_ext_sec}

To produce a reduced basis matrix $\reducedBasis$, we select fine solution candidates $\solution{\param}$ that are important columns of the solution matrix $\solutions$ and construct $\reducedBasis$ via an SVD of those important columns. However, the objective is to minimize the number of full-order solves performed during this procedure --- and retrieving a column of $\solutions$ requires a full-order solve, which is expensive. To determine the important columns of $\solutions$ without incurring this cost, we note that the coarse-proxy solutions $\solutionCoarse{\param}$ can serve as a good proxy for their fine counterparts. That is, we can identify important columns of $\solutions$ by search for important solutions among their coarse proxies $\solutionCoarse{\param}$. Thus, our initial step is to compute the entire set of coarse-proxy solutions (or alternatively, an appropriately subsampled version thereof), which we write in matrix form as
\begin{equation}
  \solutionsCoarse \equiv \begin{bmatrix}
    \solutionCoarse{\param_1} & \solutionCoarse{\param_2} & \dots & \solutionCoarse{\param_\sampleSetSize}
  \end{bmatrix}.
\end{equation} 
Note that this step is embarrassingly paralellizable. Once these solutions are ready, we identify the important elements of the sample space $\paramSet$ via a column pivoted QR decomposition of $\solutionsCoarse$. This procedure returns a permutation $\permutation$ of the columns of $\solutionsCoarse$. Let the \emph{skeleton indices} $\mathfrak{S}$ be the set of columns indices in $\permutation$ whose corresponding diagonal $R_{ii}$ is less than a certain threshold $\epsilon$ of $R_{11}$. Let the parameters $\skel{\param}_i$ corresponding to these indices be denoted as the set of \emph{skeleton parameters} $\skel{\paramSet} \subset \paramSet$. These will be our approximation as to the important columns of $\solutions$.

\subsection{Skeleton Extraction Implementation} 
We give a concrete implementation of the skeleton extraction algorithm
described above in \cref{basis_ext_alg}. This method takes in a sample space $\paramSet$ and extracts the important skeleton
parameters $\hat{\param}_j$. It returns the set of indices $\mathfrak{S} = \{i_1, i_2, ..., i_s\}$
corresponding to the indices of these skeletons, i.e. $\hat{\param}_j = \param_{i_j}$. It is
possible that the implementation can be better tailored to the problem, but we provide this
algorithm as a general-purpose default.

\begin{algorithm}[ht]
\SetAlgoLined
\KwData{A sample space $\paramSet$.}
\KwResult{The indices $\mathfrak{S}$ of the important skeleton parameters.}
\tcc{Construct coarse-proxy solutions $\solutionsCoarse$.} 
\For{$\param_i$ in $\paramSet$}{ 
    $\operatorCoarseNP(\param_i) \gets \textsc{CoarseOperator}(\param_i)$\;
    $\ve{f}(\param_i) \gets \textsc{CoarseSourceTerm}(\param_i)$\;
    $\solutionsCoarse(:, i) \gets \operatorCoarseNP(\param_i)^{-1} \ve{f}(\param_i)\;$
 } 
 \tcc{Perform column pivoted QR factorization on $\solutionsCoarse$ and denote the column
 permutation of the CPQR factorization by $\rho$.}
 $(\reducedBasis_C, \mat{R}_C, \rho) \gets \textsc{CPQR}(\solutionsCoarse)$\;
  \tcc{Select all important column indices $\rho_i$ based on $\mat{R}_{C, ii}$.}
 $\mathfrak{S} \gets \rho(\{ i \mid \mat{R}_{C, ii} \geq \epsilon \mat{R}_{C, 11} \})$\;
\Return $\mathfrak{S}$\;
\caption{\label{basis_ext_alg}$\textsc{GetSkeletons}$: Skeleton Extraction with a Coarse-Proxy Model (Offline).}
\end{algorithm}

\subsection{Reduced Basis Construction} \label{rb_construct}
Once we have selected the skeletons $\hat{\paramSet}$, we calculate the corresponding solutions for the
full-order model. We denote these corresponding fine solutions denoted as the \emph{fine skeleton set},
\begin{equation} \label{fine_solutions}
  \skel{\solutions} \equiv \solutions(:, \skel{\paramSet}) = 
  \begin{bmatrix}
    \solution{\skel{\param}_1} & \solution{\skel{\param}_2} & \dots &
    \solution{\skel{\param}_\skelCount})
  \end{bmatrix}.
\end{equation} 
Note that this step is once again embarrassingly parallelizable. 

To compute the reduced basis, we apply an SVD decomposition to the fine skeleton set
$\skel{\solutions}$ to obtain $\mat{U} \mat{\Sigma} \mat{V}^\T = \skel{\solutions}$. To build a reduced basis, we crop $\mat{U}$ by discarding all columns with singular 
values $\sigma_i$ less than $\epsilon \sigma_1$,
\begin{equation}
  \reducedBasis \equiv \mat{U}(1:\reducedBasisDim, :),
    \quad \text{ with } \reducedBasisDim \text{ such that }
  \sigma_{\reducedBasisDim} \geq \epsilon \sigma_1 > \sigma_{\reducedBasisDim + 1}\,,
\end{equation}
where $\epsilon$ is the same $\epsilon$ used in \cref{basis_ext_alg}.

We take a moment to note that the coarse-proxy model is only used to select the skeleton parameters and not for the actual construction of
the reduced basis. Hence, it is sufficient for the coarse-proxy model to be good enough to
capture the main features of the full-order model and for the important columns of the coarse-proxy solution matrix $\solutionsCoarse$ to roughly correspond to the important columns of the 
fine solution matrix $\solutions$.

\subsection{Reduced Operator Construction} \label{redops}

Once we have constructed the reduced basis $\reducedBasis$, it remains to solve the projected problem
\begin{equation}
      [ \reducedBasis^\T \operator{\param} \reducedBasis] \, \rbCoeffParam{\param} = [\reducedBasis^\T \source] \,,
  \quad
  \solution{\param} \approx \reducedBasis \, \rbCoeffParam{\param}
\end{equation}
for arbitrary $\param$. As such, we require a fast method of assembling the projected operator
$\reducedBasis^\T \operator{\param} \reducedBasis$. Assembling the full operator $\operator{\param}$
and then projecting it is prohibitively expensive. However, in solving for the fine solutions
$\skel{\solutions}$ in \cref{fine_solutions}, we have already assembled a subset of the operators
$\operator{\param}$. As we will see, the operators assembled during \cref{fine_solutions} can be
used to construct arbitrary $\reducedBasis \operator{\param} \reducedBasis^\T$ via interpolation.

If the matrix of the vectorized fine operators $\operator{\param}$ is denoted by
\begin{equation}
  \mathcal{L} \equiv \begin{bmatrix}
    \text{vec}(\operator{\param_1}) & \text{vec}(\operator{\param_2}) &
    \dots & \text{vec}(\operator{\param_\sampleSetSize})
  \end{bmatrix},
\end{equation}
then, in the process of computing the fine skeleton set $\skel{\solutions}$, we have already
assembled a subset of the columns of $\mathcal{L}$, given by
\begin{equation} \label{naivesk}
  \skel{\mathcal{L}} \equiv \mathcal{L}(:, \skel{\paramSet}) = \begin{bmatrix}
    \text{vec}(\operator{\skel{\param}_1}) & \text{vec}(\operator{\skel{\param}_2}) & \dots &
    \text{vec}(\operator{\skel{\param}_{\skelCount}})
  \end{bmatrix}.
\end{equation}
Since the operators $\operator{\param}$ change only slightly with the parameter $\param$, it stands to reason that it should be possible to use the operators we've already constructed to somehow assemble arbitrary columns of the full set of projected operators,
\begin{equation}
  \mathcal{L}_{rb} \equiv \begin{bmatrix}
    \text{vec}(\reducedBasis^\T \operator{\param_1} \reducedBasis) & 
    \text{vec}(\reducedBasis^\T \operator{\param_2} \reducedBasis) & \dots &
    \text{vec}(\reducedBasis^\T \operator{\param_\sampleSetSize} \reducedBasis) 
  \end{bmatrix}.
\end{equation}
We propose a linear interpolation method based on random samples of the fine operators $\operator{\param}$. 

To motivate our method, we first make an affine assumption. That is, we assume it is possible to assemble the $\operator{\param}$ by interpolating between the skeleton operators $\operator{\skel{\param}}$ in $\skel{\mathcal{L}}$ as such,
\begin{equation} \label{operator_interp}
    \operator{\param_i} \approx \sum_{j=1}^\skelCount \operator{\skel{\param}_j} \, \mixingCoeff_{j i}.
\end{equation}
It follows by linearity, that
\begin{equation} \label{operator_interp_proj}
    \reducedBasis^\T \operator{\param_i} \reducedBasis \approx \sum_{j=1}^\skelCount 
    \reducedBasis^\T \operator{\skel{\param}_j} \reducedBasis \, \mixingCoeff_{ji}.
\end{equation}
Note that both \eqref{operator_interp} and \eqref{operator_interp_proj} can be written in matrix form,
\begin{align} \label{least_squares}
    \mathcal{L} &\approx \skel{\mathcal{L}} \mixing, \\ 
    \mathcal{L}_{rb} &\approx \skel{\mathcal{L}}_{rb} \mixing,
\end{align}
where $\skel{\mathcal{L}}_{rb} \equiv \mathcal{L}_{rb}(:, \mathfrak{S})$ are the fine skeleton
operators projected into the reduced basis space and the $\mixing =
(\mixingCoeff_{ji}) \in \mathbb{R}^{\skelCount \times \sampleSetSize}$ is the \emph{mixing matrix} of interpolation coefficients. 
However, we must now consider how to actually compute such a mixing matrix $\mixing$. 

Our answer is based on the observation that if one makes the affine assumption, then to recover the
coefficients $\mixingCoeff_{ji}$, it suffices to randomly subsample important parts of the operators
(i.e., rows of $\mathcal{L}$) and use the resulting samples to perform least squares
regression to obtain $\mixing$. Let these important samples / row indices be denoted by
$\mathfrak{O}$. The choice $\mathfrak{O}$ can be heavily dependent on the application. For example,
if our operators are diagonally dominant, then it would make sense to include the diagonal of the
fine operators in $\mathfrak{O}$. We can also select $\mathfrak{O}$ to be slices of the operator,
which are cheap to construct, like a randomly chosen set of columns in the fine operators. Ideally,
we should have $|\mathfrak{O}| \ll \reducedBasisDim^2$. 

Taking the rows corresponding to $\mathfrak{O}$ in the above \eqref{least_squares} gives
\begin{equation} \label{least_sq_eq}
    \mathcal{L}(\mathfrak{O}, :) \approx \skel{\mathcal{L}}(\mathfrak{O}, :) \mixing.
\end{equation}
After computing $\mathcal{L}(\mathfrak{O}, :)$ for all fine operators, we then construct the mixing
matrix $\mixing$ via least-squares regression on \eqref{least_sq_eq}. Once $\mixing$ is
constructed, we can assemble any projected operator $\reducedBasis \operator{\param} \reducedBasis^\T$
by performing the linear interpolation given by \eqref{operator_interp_proj}.

\subsection{Reduced Basis and Mixing Matrix Construction Implementation}
Here, we provide an example implementation of both the construction of
the reduced basis, as described in \cref{rb_construct}, as well as the construction of the
mixing matrix described in \cref{redops}. The pseudo-code for this example implementation is given in  \cref{mixing_matrix_alg}. To use the algorithm, we require that the user
implement the following primitives:
\begin{itemize}
  \item $\textsc{FineSolve}(\param_i)$: This method takes in the parameter $\param_i$ and outputs
    the corresponding fine solution $\solution{\param_i}$ as well as the corresponding vectorized fine
    operator $\text{vec}(\operator{\param_i})$. \textbf{Nota bene} that, in practice, it may be the case that $\text{vec}(\operator{\param_i})$ may be too large to store in memory. This is not an obstacle, as we only use $\text{vec}(\operator{\param_i})$ for notational convenience. To implement the following algorithms, one only needs to be able to apply the operator $\operator{\param_i}$ and to be able to sample a sparse subset of the entries of $\text{vec}(\operator{\param_i})$. 
              \item $\textsc{GetOperatorSamples}()$: This method chooses the set of operator entries
    $\mathfrak{O}$ (i.e., rows of the matrix $\mathcal{L}$) to sample and outputs the operator
    samples $\mathcal{L}(\mathfrak{O}, :)$, as described in \cref{redops}.
\end{itemize}
Note that there is a part of the implementation which involves adding additional skeletons to the
skeleton set. This segment of the algorithm will be addressed in \cref{add_sk}.

\begin{algorithm}[ht]
\KwData{A sample space $\paramSet$}
\KwResult{A reduced basis matrix $\reducedBasis$, a mixing matrix $\mixing$, projections
$\skel{\mathcal{L}}_{rb}$ of fine skeleton operators into the reduced basis space.}
\tcc{Compute the important skeletons in the sample space}
$\mathfrak{S} \gets \textsc{GetSkeletons}(\paramSet)$\;
$\skel{\paramSet} \gets \paramSet(:, \mathfrak{S})$\;
\tcc{Compute the corresponding fine skeleton solutions}
\For{$\skel{\param}_j$ in $\skel{\paramSet}$}{
  $(\skel{\solutions}(:, j), \skel{\mathcal{L}}(:, j)) \gets \textsc{FineSolve}(\skel{\param}_j)$\;
}
\tcc{(Optional) Use additional skeleton extraction as described in \cref{add_sk}}
\If{Using additional skeleton extraction}{
  $(\skel{\solutions}(:, j), \skel{\mathcal{L}}(:, j)) \gets 
  \textsc{AdditionalSkeletons}(\paramSet, \skel{\solutions}, \skel{\mathcal{L}}, \mathfrak{S}, \mathcal{L}_{samp})$\;
}

\tcc{Construct reduced basis matrix $\reducedBasis$ from fine skeletons $\skel{\solutions}$ by taking the first few left singular vectors of $\skel{\solutions}$.}
$(\mat{U}, \mat{\Sigma}, \mat{V}) \gets \textsc{SVD}(\skel{\solutions})$\;
$\reducedBasis \gets \mat{U}(:, \mat{\Sigma} > \epsilon \sigma_1)$\;
\tcc{Compute the samples $\mathcal{L}(\mathfrak{O}, :)$ from each fine operators}
$\mathcal{L}_{samp} \gets \textsc{GetOperatorSamples}()$\;
\tcc{Perform least squares regression using the samples $\mathcal{L}_{samp}$ to compute the mixing matrix $\mixing$.}
$\mixing \gets \textsc{LeastSquares}(\mathcal{L}_{samp}, \mathcal{L}_{samp}(:, \mathfrak{S}))$\;
\tcc{Project the skeleton operators $\skel{\mathcal{L}}$ into the reduced basis space given by $\reducedBasis$. Note $L(\skel{\param}_j)$ has been reshaped into a matrix.}
\For{$\text{vec}(\operator{\skel{\param}_j})$ in $\skel{\mathcal{L}}$}{
  $\skel{\mathcal{L}}_{rb}(:, j) \gets \text{vec}(\reducedBasis^\T \operator{\skel{\param}_j} \reducedBasis)$\;
}
\Return $(\reducedBasis, \mixing, \skel{\mathcal{L}}_{rb})$\;
\caption{\label{mixing_matrix_alg}Reduced Basis and Mixing Matrix Computation. (Offline)}
\end{algorithm}

\subsection{Online Reduced Basis Solve Implementation}
In this subsection, we provide pseudo-code in \cref{reduced_basis_alg} for using the offline
computations performed in \cref{basis_ext_alg,mixing_matrix_alg} to compute fast online
approximations to $\solution{\param}$ for arbitrary $\param \in \paramSet$. We suppose that we are
provided with the following primitive:

\begin{algorithm}[ht] 
\KwData{A sample $\param \in \paramSet$ for which to compute a reduced basis approximation, the mixing matrix $\mixing$, the projected skeleton operators $\skel{\mathcal{L}}_{rb}$, and the reduced basis matrix $\reducedBasis$.}
\KwResult{An approximation $\rbApprox$ of $\solution{\param}$.}
\tcc{Assemble our approximation for the projected operator $\operatorRBNP \equiv \reducedBasis^\T \operator{\param} \reducedBasis$ using the projected skeleton operators $\skel{\mathcal{L}}_{rb}$ and the mixing matrix $\mixing$.}
$\text{vec}(\operatorRBNP) \gets \skel{\mathcal{L}}_{rb} \, \mixing(:, i)$\;
\tcc{Have the oracle assemble the right hand side of the equation, i.e., $\sourceRB(\param) \equiv
\reducedBasis^\T \source$, for us and project it into the reduced basis space.}
$\sourceRB \gets \textsc{AssembleRightHandSide}(\reducedBasis, \param)$\;
\tcc{Solve the system and return the result.}
$\rbCoeff \gets \operatorRBNP^{-1} \sourceRB$\;
\tcc{Lift result from reduced basis space to $\mathbb{R}^n$.}
$\rbApprox \gets \reducedBasis \, \rbCoeff $\;
\Return $\rbApprox$\;
\caption{\label{reduced_basis_alg}Reduced Basis Solve for $\solution{\param}$ (Online).}
\end{algorithm}

\begin{itemize}
  \item $\textsc{AssembleRightHandSide}(\reducedBasis, \param_i)$: This method takes in the reduced
    basis $\reducedBasis$ and a paramter $\param_i$ and returns $\reducedBasis^\T \ve{f}(\param_i)$ or an
    approximation thereof. Depending on the problem being solved, there might be some intricacies to
    this. However, if $\ve{f}(\param_i)$ is inexpensive to assemble, then the oracle can simply compute
    $\ve{f}(\param_i)$ and apply $\reducedBasis^\T$. In other situations, one can use mathematical
    manipulations to obtain an expression for $\ve{f}(\param_i)$ in terms of already computed expressions.
    See the radiative transport equation \cref{retsec} for a nontrivial case. In the worst
    case, if the entries of $\ve{f}(\param_i)$ are not overly expensive to sample, one can sub-sample the
    $\ve{f}(\param_i)$ and use the samples to linearly interpolate between $\reducedBasis^\T
    \ve{f}(\skel{\param}_i)$ by constructing a mixing matrix using the technique in \cref{redops}. A
    more involved sub-sampling alternative could be to use a discrete empirical interpolation method
    such as Q-DEIM \cite{drmac2016new} to compute sub-sampling entries in $\ve{f}(\omega)$ and
    interpolation weights for $\reducedBasis^\T \ve{f}(\omega)$.
\end{itemize}

\subsection{A Note on Gappy Matrix POD} We remark that the above method of constructing reduced
operators is close to Gappy Matrix POD in \cite{carlberg2015preserving}. However, one key
distinction is that we do not orthogonalize the skeleton operators $\reducedBasis^\T
\operator{\skel{\param}_i} \reducedBasis$. Gappy Matrix POD would involve vectorizing the skeleton operators
$\operator{\skel{\param}_i}$, taking SVD to find a set of orthogonalized operators
$\operatorNP^{\perp}_1, ..., \operatorNP^{\perp}_r$, projecting them into the reduced basis space, and
then using $\reducedBasis^\T \operatorNP^{\perp}_1 \reducedBasis, ... \reducedBasis^\T
\operatorNP^{\perp}_r \reducedBasis$ to interpolate the general projected operators
$\reducedBasis^\T \operator{\param_i} \reducedBasis$. We do not do this. This is intentional. While performing this
orthogonalization may sometimes result in increased stability of interpolation, for integral
operators, it is not desirable to represent the underlying operators as full dense matrices.
Moreover, by virtue of how we select the skeleton operators, we ensure to some extent
that the interpolation problem is already relatively well-conditioned.

\subsection{Additional Skeleton Extraction} \label{add_sk}

Sometimes, the fine operators $\skel{\mathcal{L}}$ we assemble during our fine solves may not be sufficiently rich to reconstruct all of the operators in $\mathcal{L}_{rb}$ via interpolation. If this is the case, then we must add additional columns to our set of skeleton operators $\skel{\mathcal{L}}$. Note that we can use the fine operator samples $\mathcal{L}(\mathfrak{O}, :)$ in the previous section to get a rough idea the important operators in $\mathcal{L}$. To find the operators we have failed to represent well with our choice of skeletons $\mathcal{L}$, we can consider the operator samples $\mathcal{L}(\mathfrak{O}, :)$ with our skeletons $\mathcal{L}(\mathfrak{O}, \mathfrak{S})$ projected out,
\begin{equation}
    \mathcal{L}_{res} \equiv \mathcal{L}(\mathfrak{O}, :) - P \mathcal{L}(\mathfrak{O}, :) ,
\end{equation}
where $P$ is a projector onto the column space of $\mathcal{L}(\mathfrak{O}, \mathfrak{S})$. We call these the \emph{residual operator samples}. This projection can be done via modified Gramm-Schmidt, for example. 

Then, before we compute the mixing matrix, we can perform a column pivoted QR decomposition of
$\mathcal{L}_{res}$ to find operators we're unable to approximate well. Similar to what was done in
\cref{rb_construct}, we select the columns with a diagonal $R$-factor which is smaller than $\epsilonMul \epsilon$ multiplied by the largest column norm in the unprojected $\mathcal{L}(\mathfrak{O}, :)$, where $\epsilonMul$ is an arbitrary constant set by the user. Whatever columns $\mathfrak{A}$ are selected by this process, we append them to our set of fine operator skeletons $\skel{\mathcal{L}}$ as such,
\begin{equation} \label{appendsk}
  \skel{\mathcal{L}} \gets \begin{bmatrix}
    \skel{\mathcal{L}} & \mathcal{L}(:, \mathfrak{A})
  \end{bmatrix}.
\end{equation}
In addition, depending on the problem at hand, one can also add the corresponding fine solutions of $\mathfrak{A}$ to the fine solution skeleton set $\skel{\solutions}$, as these may add important fine scale information which our coarse-proxy model may have missed,
\begin{equation} \label{appendsl}
  \skel{\solutions} \gets \begin{bmatrix}
    \skel{\solutions} & \solutions(:, \mathfrak{A})
  \end{bmatrix}.
\end{equation}
Afterwards, one can continue with everything detailed in \cref{redops} without any changes, using
\cref{appendsk} instead of \cref{naivesk} for the skeleton operators $\mathcal{L}$. 

\subsection{Implementation of Additional Skeleton Extraction}

We now provide a pseudo-code implementation in \cref{add_sk_alg} of the additional skeleton
extraction algorithm presented above in \cref{add_sk}.

\begin{algorithm}[ht]
\SetAlgoLined
\KwData{A sample space $\paramSet$, a set of fine skeletons $\skel{\solutions}$, their corresponding operator skeletons $\skel{\mathcal{L}}$ and indices $\mathfrak{S}$, and a matrix of operator samples $\mathcal{L}_{samp}$.}
\KwResult{A possibly enlarged set of operator skeletons $\skel{\mathcal{L}}$ and fine solution skeletons $\skel{\solutions}$.}
\tcc{Compute the maximum energy in the operator samples before we project out the fine skeletons.}
$a \gets \max_i \| \mathcal{L}_{samp}(:, i) \|_2$\;
\tcc{Project out the fine skeletons $\mathcal{L}_{samp}(:, \mathfrak{S})$ we've computed in algorithm 1 from the samples $\mathcal{L}_{samp}$.}
$\mathcal{L}_{res} \gets \textsc{ProjectOut}(\mathcal{L}_{samp}, \mathcal{L}_{samp}(:, \mathfrak{S}))$\;
\tcc{Extract important operators we've missed during reduced basis extraction via QR decomposition of residual samples $\mathcal{L}_{res}$.}
$(\reducedBasis, R, \rho) \gets \textsc{CPQR}(\mathcal{L}_{res})$\;
\tcc{Select only column indices for which $r_{ii} \geq \epsilonMul \epsilon a$. In practice, this should be done by stopping the above QR factorization when this first happens.}
$\mathfrak{A} \gets \rho(r_{ii} \geq \epsilonMul\epsilon a)$\;

\tcc{ Compute new set of additional fine operators $\skel{\mathcal{L}}_A$ and additional fine solutions $\skel{\solutions}_A$ for the selected columns $\mathfrak{A}$.}
\For{$\param_j'$ in $\paramSet(:, \mathfrak{A})$}{
$(\skel{\solutions}_A(:, j), \skel{\mathcal{L}}_A(:, j)) \gets \textsc{FineSolve}(\param_j')$\;
}
\tcc{Add additional skeleton solutions $\skel{\solutions}_A$ to our existing skeleton solutions $\skel{\solutions}$.}
$\skel{\solutions} \gets \begin{bmatrix} \skel{\solutions} & \skel{\solutions}_A \end{bmatrix}$\;
\tcc{Append new skeletons $ \skel{\mathcal{L}}_A$ to our existing skeletons $\skel{\mathcal{L}}$}
$\skel{\mathcal{L}} \gets \begin{bmatrix} \skel{\mathcal{L}} & \skel{\mathcal{L}}_A \end{bmatrix}$\;
\Return $(\skel{\solutions}, \skel{\mathcal{L}})$\;
\caption{\label{add_sk_alg}$\textsc{AdditionalSkeletons}$: Optional Additional Skeleton Extraction (Offline). }
\end{algorithm}

\subsection{Interpolating Operators with an Offset} \label{offset_ops}

There are many problems in which the operators $\operator{\param}$ take on a natural form,
\begin{equation}
    \operator{\param} = \operatorAffineOffset + \operatorAffineParam{\param} ,
\end{equation}
where $\operatorAffineOffset$ does not depend explicitly on $\param$, and is shared among all of the
operators $\operator{\param}$. We will see such examples of this later. In such situations, it may
be more advisable to interpolate the operator $\operatorAffineParam{\param}$ instead of the full
operator $\operator{\param}$ when constructing reduced operators. All techniques from \cref{redops} carry over with minimal modification. One assumes that there exist interpolation coefficients for the operators $\operator{\param}$,
\begin{equation} \label{operator_interp2}
    \operatorAffineParam{\param_i} \approx \sum_{j} \alpha_{i j} \operatorAffineParam{\skel{\param}_j} .
\end{equation}
Then one can find a mixing matrix $\mixing$ with
\begin{equation}\label{least_squares2}
    \mathcal{B}_{rb} \approx \skel{\mathcal{B}}_{rb} \mixing ,
\end{equation}
by simply performing a least squares solve of the equation
\begin{equation} \label{least_sq_eq2}
    \mathcal{B}(\mathfrak{O}, :) \approx \skel{\mathcal{B}}(\mathfrak{O}, :) \mixing .
\end{equation}
where $\mathcal{B}$, $\skel{\mathcal{B}}$, $\mathcal{B}_{rb}$, $\skel{\mathcal{B}}_{rb}$ are defined
analogously to $\mathcal{L}$, $\skel{\mathcal{L}}$, $\mathcal{L}_{rb}$, and
$\skel{\mathcal{L}}_{rb}$ in \cref{redops}. 

Once the mixing matrix $\mixing$ has been computed, note that the corresponding coefficients
$\mixingCoeff_{ji}$ can then be used to interpolate the reduced operators $\reducedBasis^\T \operator{\param_i} \reducedBasis$,
\begin{equation}
  \reducedBasis^\T \operator{\param_i} \reducedBasis \approx \reducedBasis^\T
  \operatorAffineOffset \reducedBasis + 
  \sum_{j=1}^\skelCount \reducedBasis^\T \operatorAffineParam{\skel{\param}_j} \reducedBasis \mixingCoeff_{ji} .
\end{equation}
The quantity $\reducedBasis^\T \operatorAffineOffset \reducedBasis$ can be computed alongside the
skeleton operators $\reducedBasis^\T \operatorAffineParam{\skel{\param}_j} \reducedBasis$ when
projecting operators into the reduced basis space.

\section{Numerical results}\label{sec:num}
To demonstrate that our framework is both a practical and efficient approach to model order
reduction for integral equations, we perform simulations on the two following examples.

\subsection{Boundary Integral Formulation of the Laplace Equation}
We consider the Laplace equation
\begin{equation}
  \begin{aligned}
    \Delta \varphi &= 0, \; \text{ in } \domain, \\
    \varphi &= f, \; \text{ on } \partial \domain,
  \end{aligned}
\end{equation}
where $\domain\subset \mathbb{R}^2$ is a bounded Lipschitz domain. Introduce the single layer
potential $u$, which is given by the solution to the integral equation
\begin{equation} \label{bie}
    f(\bm{x}) = \frac{1}{2} u(\bm{x}) - \int_{\partial \domain} \frac{\partial G(\bm{x},
    \bm{y})}{\partial n(\bm{y})} u(\bm{y})  \dd s(\bm{y}) ,
\end{equation}
where $G$ denotes the Green's function of the Laplace equation in 2D,
\begin{equation}
    G(\bm{x}, \bm{y}) =  \frac{1}{2\pi} \ln \frac{1}{|\bm{x} - \bm{y}|}.
\end{equation}
Then once $u$ has been computed by solving the integral equation \cref{bie}, the solution $\phi$ to
Laplace equation can be recovered by
\begin{equation}
    \varphi(\bm{x}) = - \int_{\partial \domain} \frac{\partial G(\bm{x}, \bm{y})}{\partial n(\bm{y})} 
    u(\bm{y}) \dd s(\bm{y}).
\end{equation}
Hence, the key of this problem is to solve the boundary integral equation \cref{bie}.

To bring this problem in the many-query setting of our framework, we suppose that the shape of the
domain $\domain(\param)$ is parameterized by $\param$ taken from a sample space $\paramSet_\infty$.
Thus, the integral equations we would like to solve are given by
\begin{equation} \label{bie_want_to_solve}
    f(\bm{x}) = \frac{1}{2} u(\bm{x}; \param) - \int_{\partial \domain(\param)} \frac{\partial
    G(\bm{x}, \bm{y})}{\partial n(\bm{y})} u(\bm{y}; \param) \dd s(\bm{y}) ,
\end{equation}
where $f(\bm{x})$ is a function prescribed on $\mathbb{R}^2$, which we hold constant across all
problem instances in $\paramSet_\infty$. One can discretize this equation by taking a discrete
number of samples of $u$ on $\partial \domain$, and using an appropriate integral quadrature
for the integral kernel above. An example oracle for this problem is one that uses a significantly
reduced number of quadrature points on $\partial \domain$ as its coarse-proxy model. Note that using a
coarse-proxy model with low resolution becomes difficult when the source function $f$ exhibits singular
behavior near the boundary $\partial \domain$, as the high frequency content in $u$ is
difficult to resolve. Nonetheless, even if such a coarse-proxy model is used for skeleton extraction, our
results suggest that using the additional skeleton extraction techniques discussed in \cref{add_sk}
will compensate for the information which the coarse-proxy model cannot resolve, since this high-frequency
information will present itself in the operator samples.

Regardless of the exact implementation of the oracle's scheme for extracting skeletons, the
discretized equation \eqref{bie_want_to_solve} reads
\begin{equation}\label{bie_ie}
    \operator{\param} \solution{\param} = \mathsf{f}(\param),
\end{equation}
where $\solution{\param}$ is the discretized version of the double layer potential $u$, $\mathsf{f}(\param)$ is the source term $f$ sampled on $\partial \domain(\param)$. And $\operator{\param}$ is a matrix that has the form
\begin{equation}
    \operator{\param} = \frac{1}{2} \mathsf{I} - \mathsf{G}(\param),
\end{equation}
where $\mathsf{I}$ is the identity, and $\mathsf{G}(\param)$ is the discretized integral kernel in
\cref{bie_want_to_solve}. 

After skeleton extraction and construction of the reduced basis $\reducedBasis$, the
projected equations in the reduced basis are given by
\begin{equation}
    \left(\frac{1}{2} \mathsf{I} - \reducedBasis^\T \mathsf{G}(\param) \reducedBasis\right)
    \rbCoeffParam{\param}= \reducedBasis^\T \mathsf{f}(\param)
    \quad
    \solution{\param} \approx \reducedBasis \, \rbCoeffParam{\param}
\end{equation}
One can construct the reduced operators $\reducedBasis^\T \mathsf{G}(\param) \reducedBasis$ via the
techniques discussed in \cref{redops}. For our operator samples, we sample a few columns of
all operators $\mathsf{G}(\param)$. In general, $\mathsf{f}(\param)$ is inexpensive to assemble, so
one can simply construct the source term $\mathsf{f}(\param)$ and project it into the reduced basis
space by applying $\reducedBasis^\T$ during an online reduced basis solve. 

\subsubsection{Results}
We parameterize the boundary $\partial \domain(\param)$ by a polar curve $\gamma(\theta; \param) :
\mathbb{R} \times \paramSet \longrightarrow \mathbb{R}^2$, where $\theta \in [0, 2 \pi)$ and
\begin{equation}
  \gamma(\theta; \param) \equiv r(\theta; \param) \left[\cos(\theta), \sin(\theta)\right]^T,
\end{equation}
where $r(\theta; \param) : \mathbb{R} \times \paramSet \longrightarrow \mathbb{R}$ is the radial
distance of this curve from the origin at angle $\theta$. To specify the radial function, we chose a
set of interpolation points $\theta = 2 \pi k /N$ for $k \in \{0, ..., N - 1\}$ and require that
\begin{equation}
  r(2 \pi k / N; \omega) = b_k(\omega)
\end{equation}
where $b_k(\omega)$ are interpolation value for the radial function $r(\theta; \param)$ at the
interpolation points $\theta = 2 \pi k / N$. We then determine remainder of the curve $\gamma$ by
Fourier interpolation. Viewing the continuous parameter set $\paramSet_\infty$ as a probability
space, we take the radial interpolation points $b_k(\param)$ to be i.i.d. uniformly random in the
interval,
\begin{equation}
  b_k(\param) \sim \mathcal{U}([1 - \kappa, 1 + \kappa]) .
\end{equation}
Finally, to make the problem challenging, we take the source function $f(\bm{x})$ to have a singularity which can potentially be situated near the curve $\gamma$. In particular, we take
\begin{equation}
f(\bm{x}) \equiv \frac{1}{\|\bm{x} - \bm{x}_0\|_2} .
\end{equation}
To ensure there are problem instances where the curve $\gamma$ comes close to the singularity located at $\bm{y}$ we take our parameters to be
\begin{equation}\label{eq:bie_setup}
\kappa = 0.4, \qquad \bm{x}_0 = (0.6, 0), \qquad N = 8 .
\end{equation}
Finally, to extract a discrete parameter space $\paramSet \subset \paramSet_\infty$, we draw
$|\Omega| = 32768$ random samples from $\paramSet_\infty$. 

For our full-order model, we use a total of $n_f = 2048$ quadrature samples. Whereas, for our coarse-proxy model, we use a total of $n_c = 128$ quadrature samples. For selecting additional
skeletons, we use a selection threshold multiplier of $\epsilonMul = 1.5$. 

\begin{figure}[htb]
 \centering 
 \subfloat[Diagonal of $R$]{
   \includegraphics[width=0.4\textwidth]{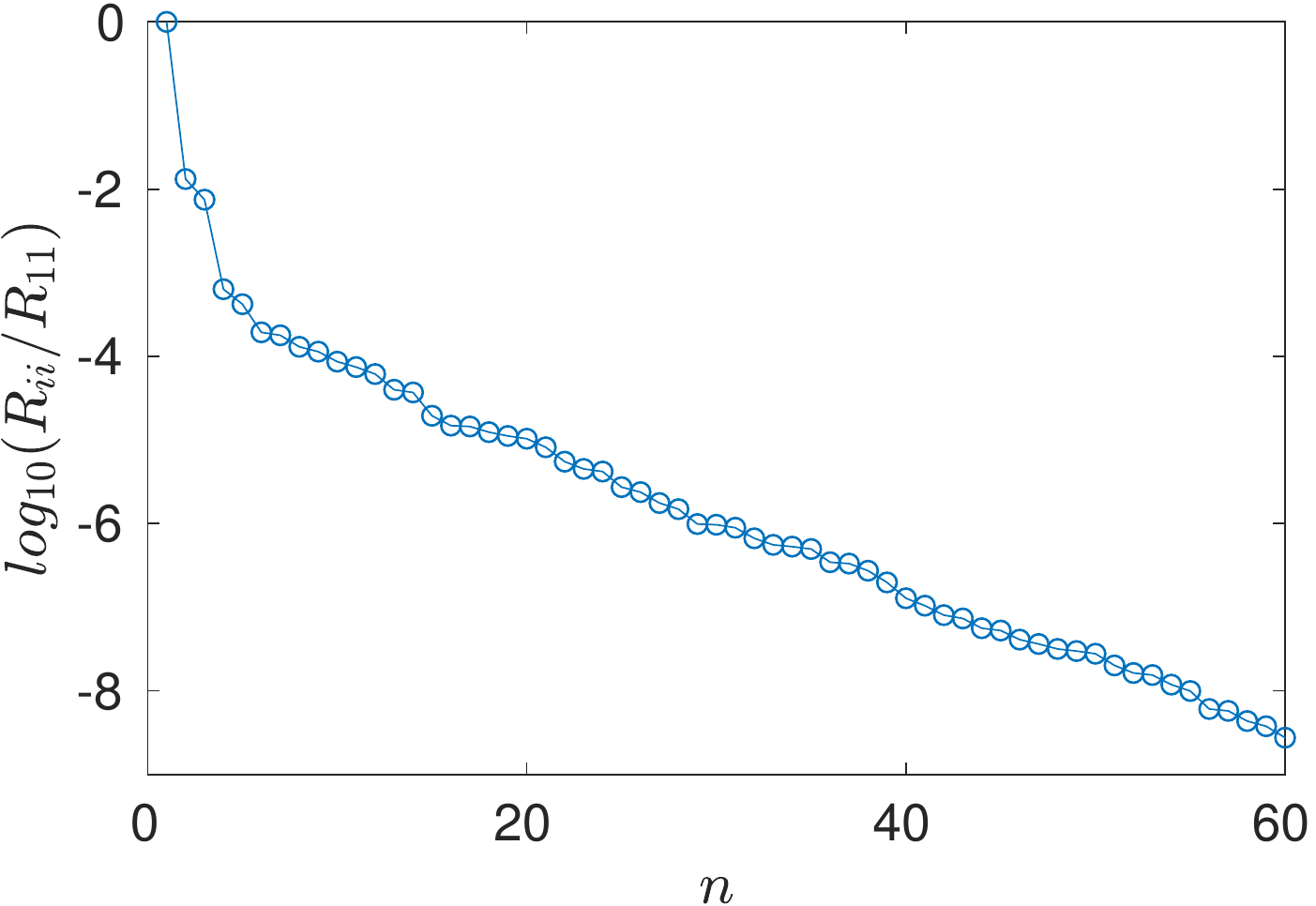}
 }\hspace{0.1\textwidth}
 \subfloat[Sigular value $\sigma_k$]{
   \includegraphics[width=0.4\textwidth]{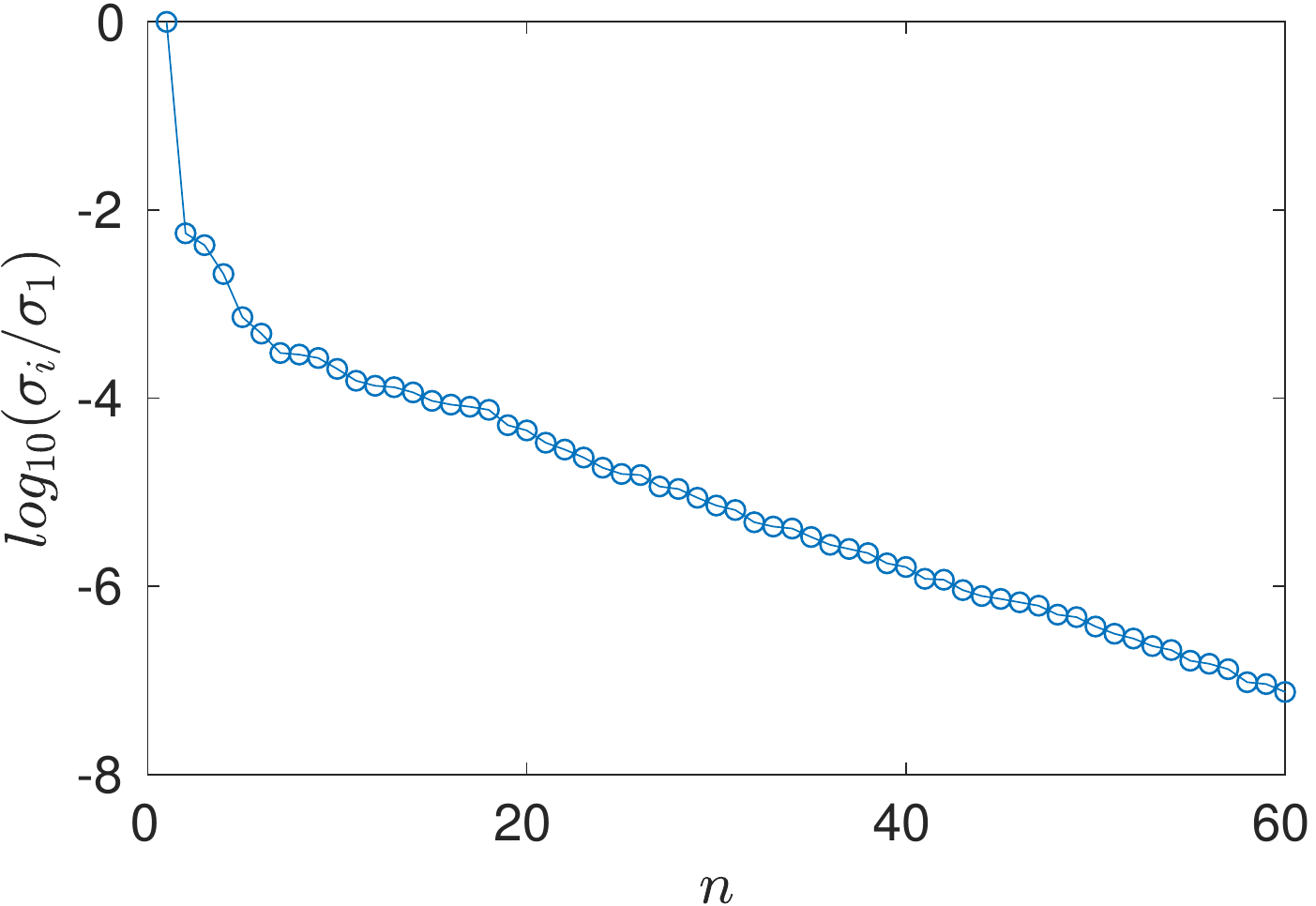}
 }
 \caption{\label{fig:bie_Rsigma} Left: the normalized diagonal values of $R$ in the skeleton selection for the described Laplace equation problem. Right:
 the normalized singular values in the basis construction for the described Laplace equation problem.}
\end{figure}

\begin{figure}[htb]
  \centering
  \includegraphics[width=0.95\textwidth]{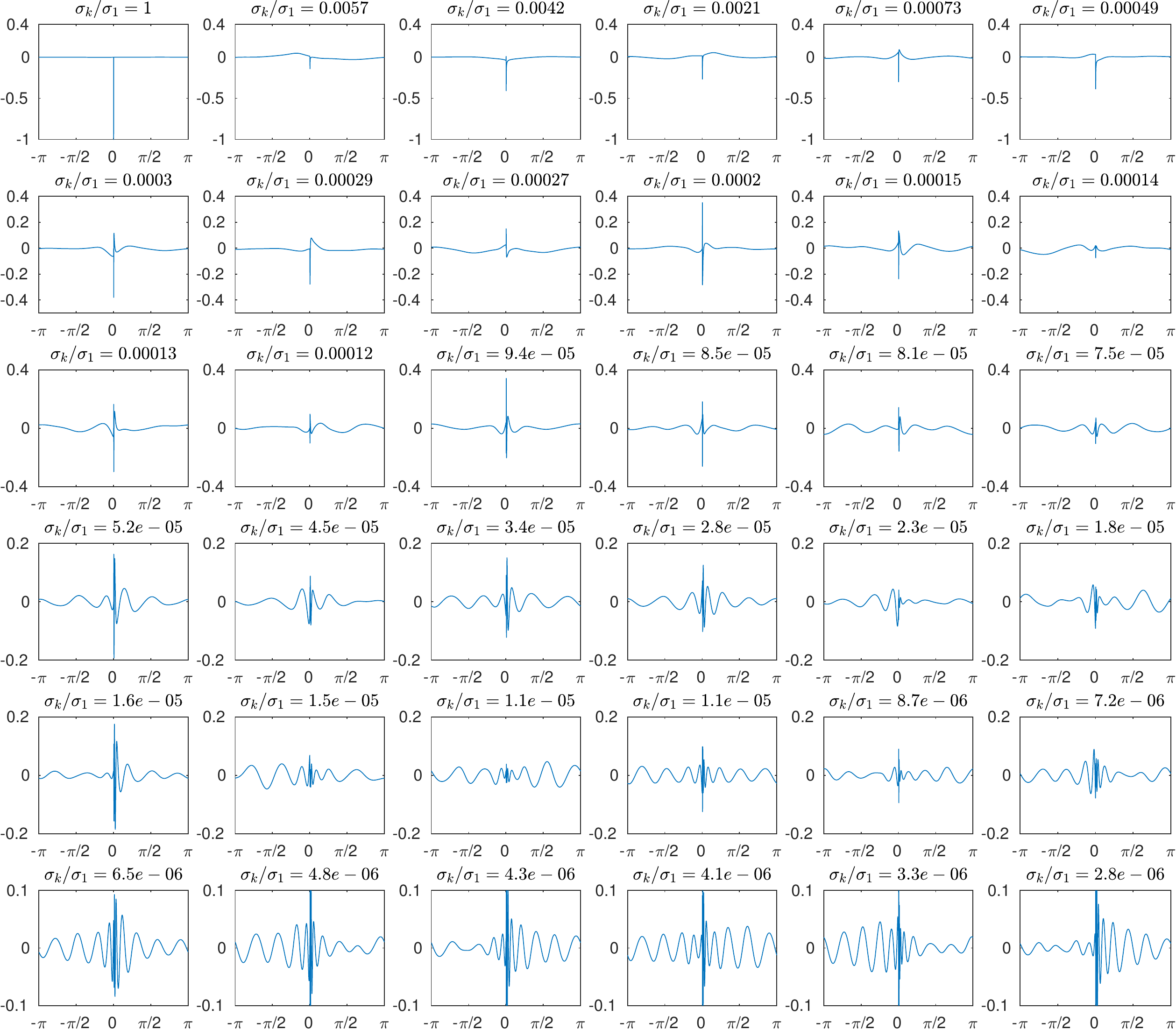}
  \caption{\label{fig:basis_bie} 
    Profiles of the reduced basis generated for the described boundary integral form of the
  Laplace equation problem for different singular values $\sigma_k$.}
\end{figure}

\begin{figure}[htb]
  \centering
  \begin{tabular}{ccc}
    \begin{overpic}[width=0.26\textwidth]{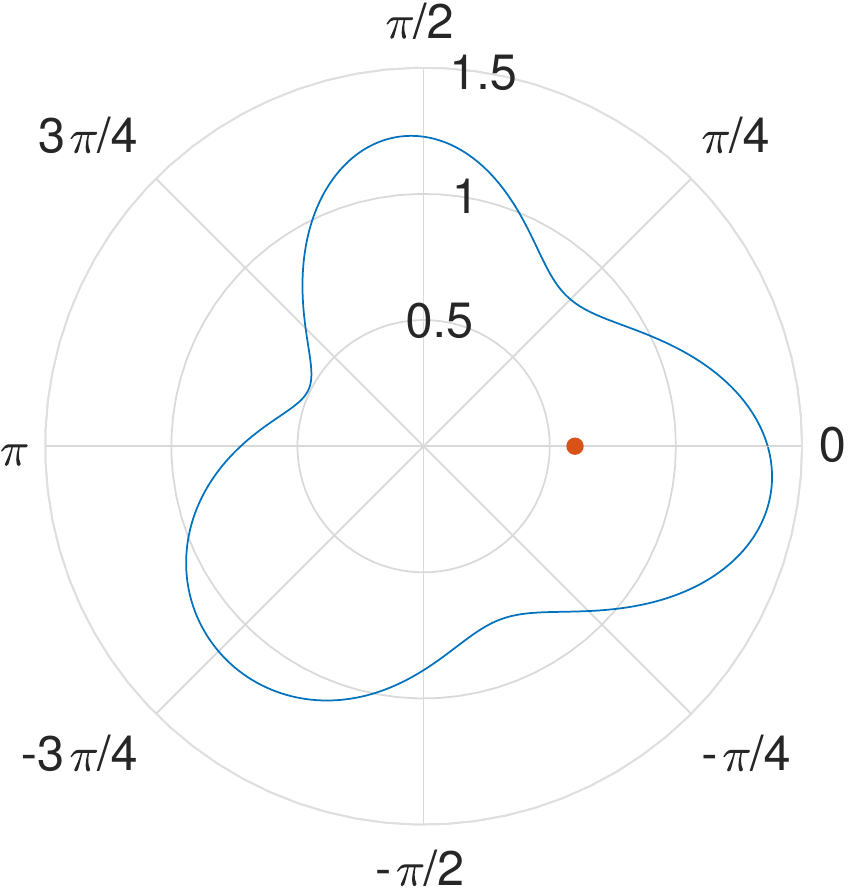} 
      \put(43,46){$\mathcal{D}$}
    \end{overpic} &
    \begin{overpic}[width=0.26\textwidth]{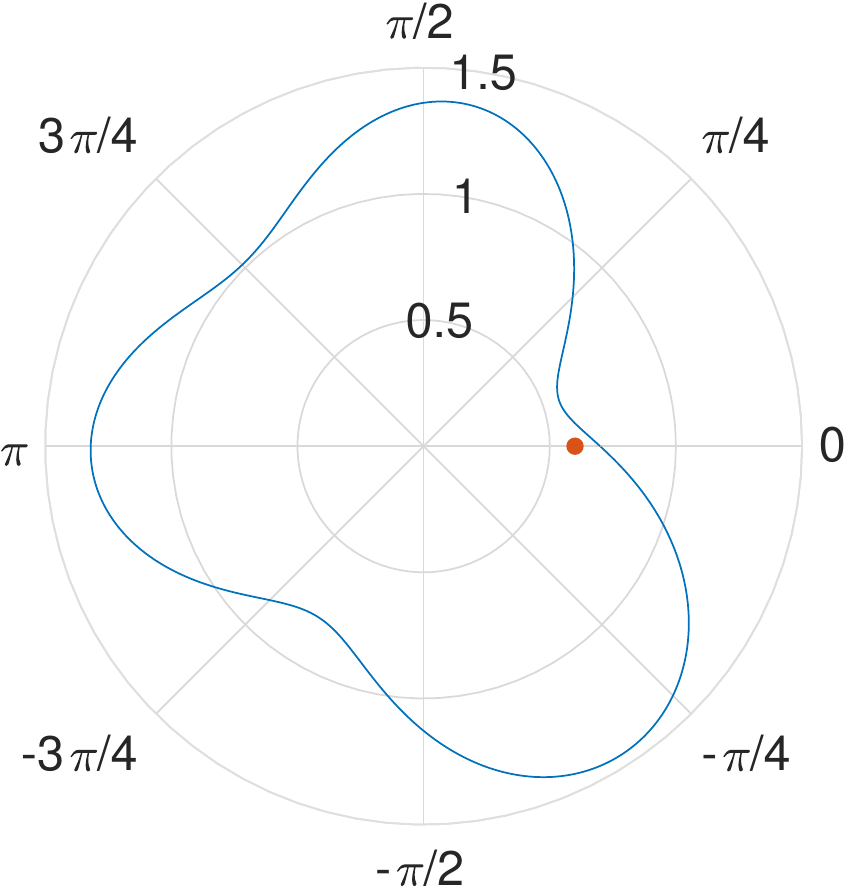} 
      \put(43,46){$\mathcal{D}$}
    \end{overpic} &
    \begin{overpic}[width=0.26\textwidth]{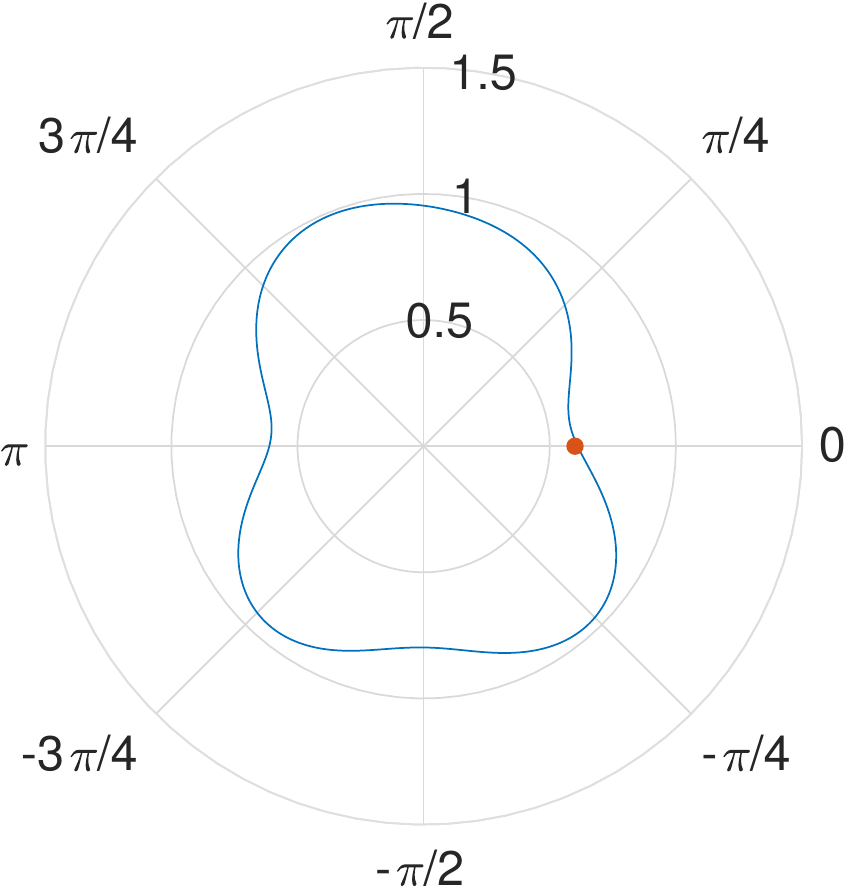} 
      \put(43,46){$\mathcal{D}$}
    \end{overpic} \\
    \includegraphics[width=0.3\textwidth]{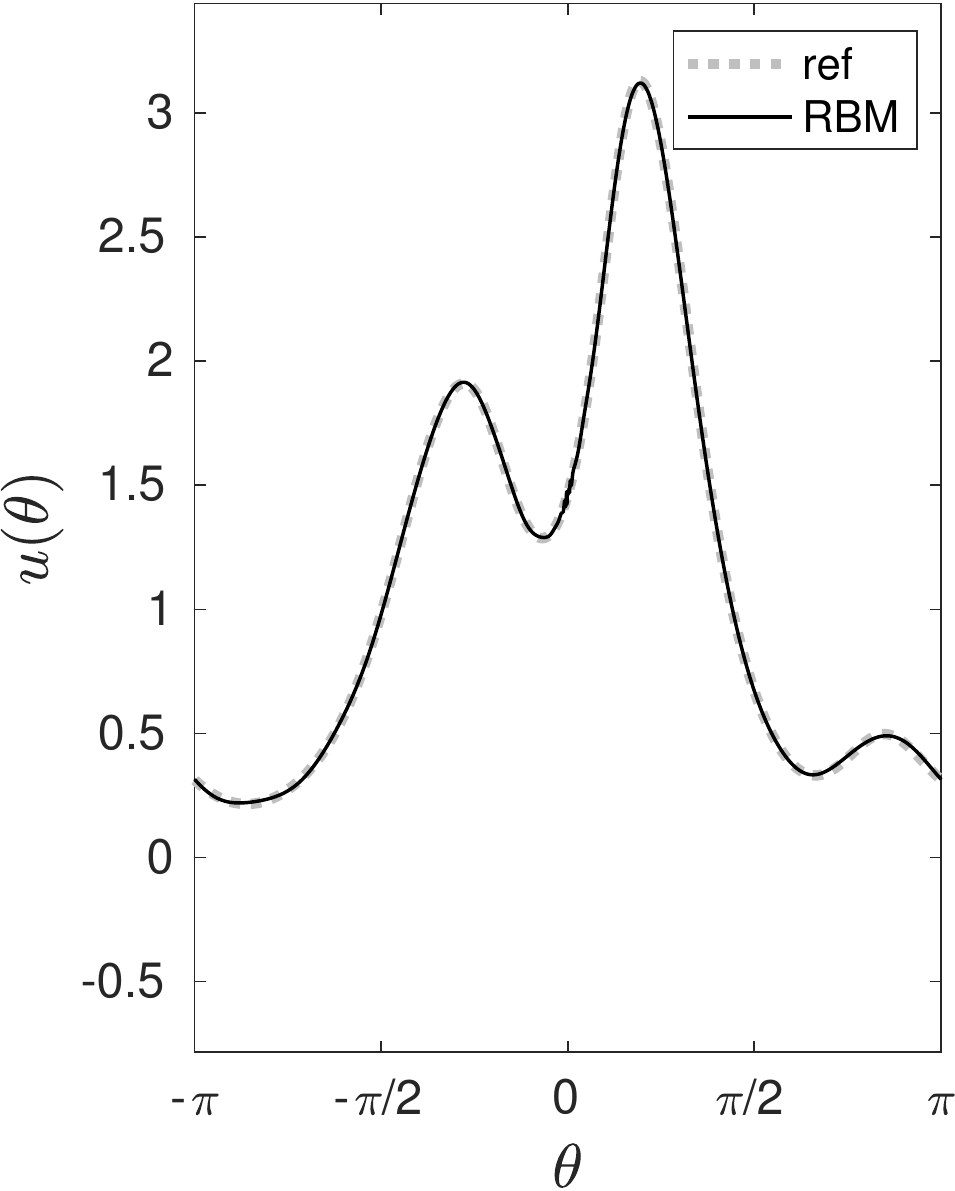} \hspace{0.02\textwidth} &
    \includegraphics[width=0.3\textwidth]{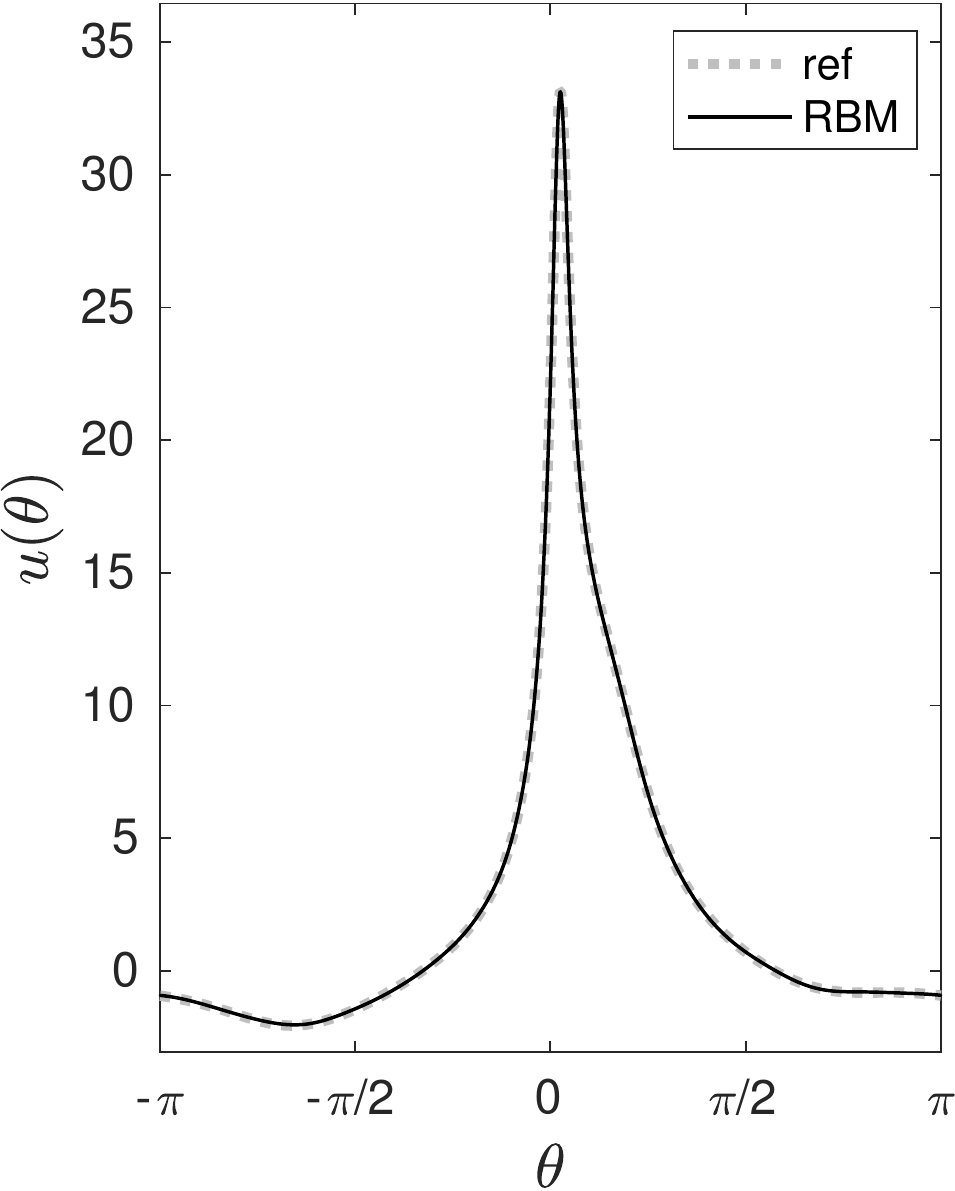} \hspace{0.02\textwidth} &
    \includegraphics[width=0.3\textwidth]{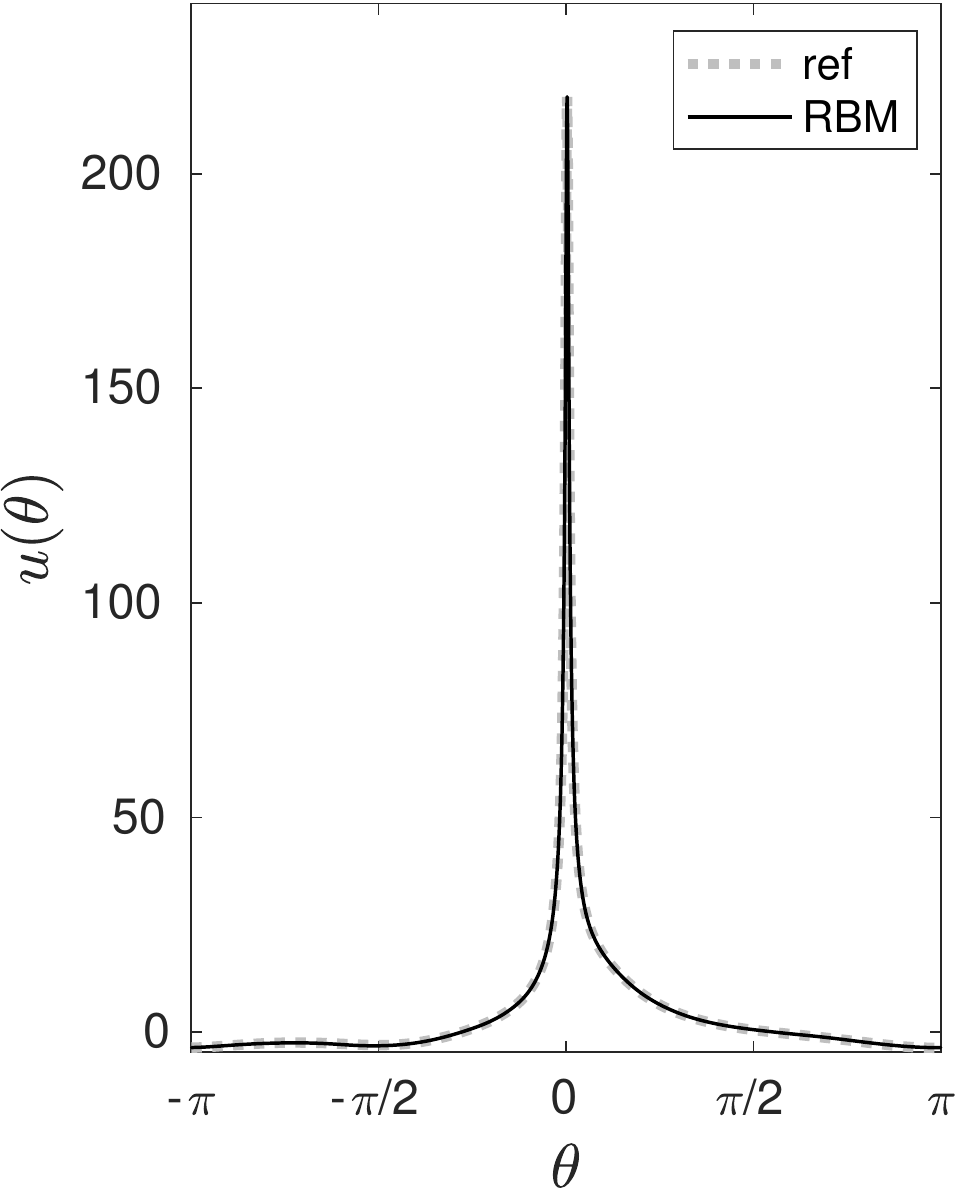} \\
    \includegraphics[width=0.3\textwidth]{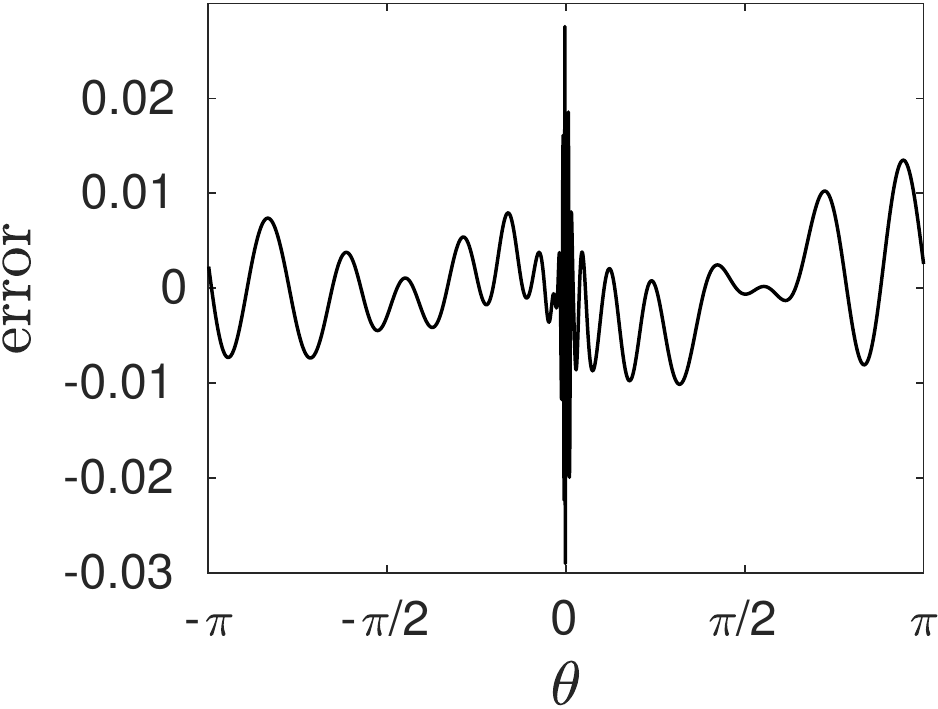} \hspace{0.02\textwidth} &
    \includegraphics[width=0.3\textwidth]{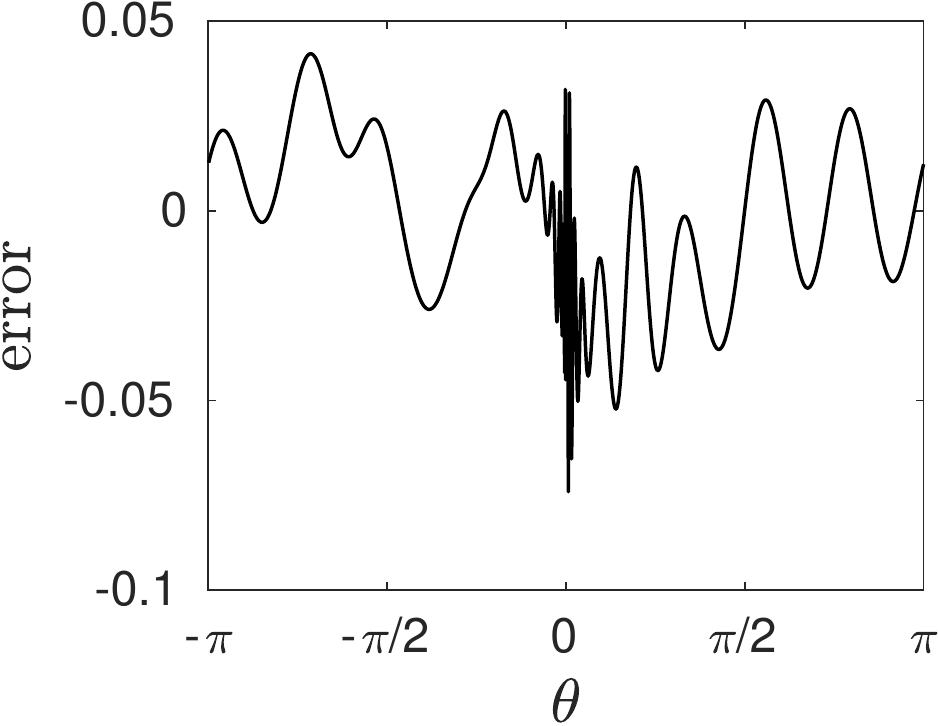} \hspace{0.02\textwidth} &
    \includegraphics[width=0.3\textwidth]{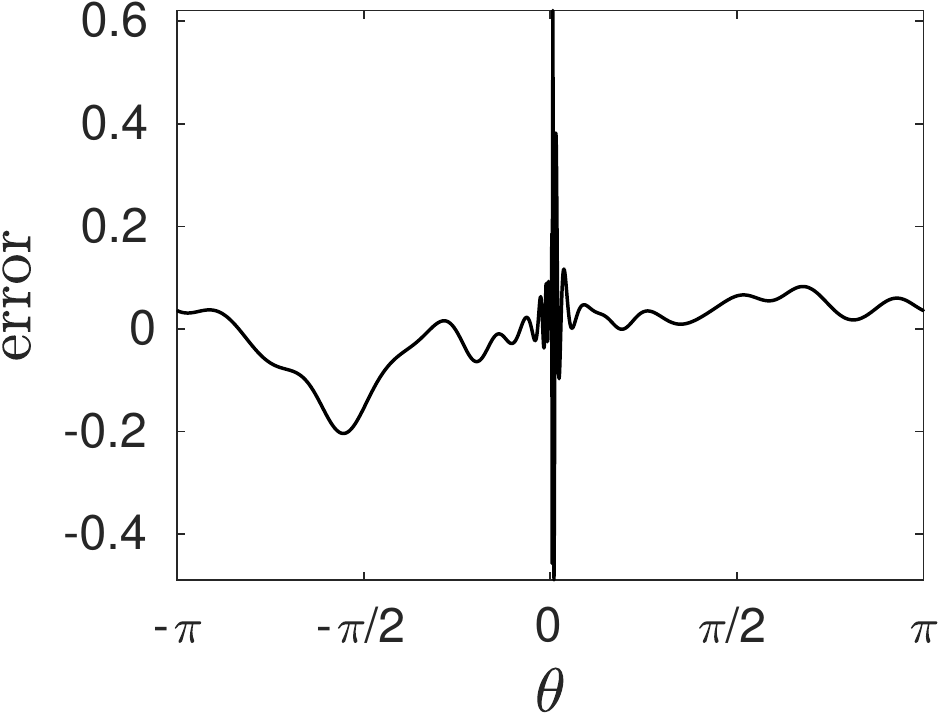} \\
  \end{tabular}
        \caption{\label{fig:bie_examples} 
    Three examples of the solutions evaluated by the reduced basis method (RBM) for the
    single layer potential $u(\theta; \omega)$ for the parameter set $\paramSet$ described for
    the boundary integral form of the Laplace equation and their corresponding reference solution
  (ref) with threshold $\epsilon=1\times 10^{-6}$. The upper figures are the domain $\mathcal{D}$
and the red point is the location of the singularity in \cref{eq:bie_setup}.}
                \end{figure}

\begin{table}[htb]
\centering
 \begin{tabular}{|| c || c | c || c | c| c|c || c ||} 
 \hline
 $\epsilon$ & $\skelCount$ & $\reducedBasisDim$ & $T_{rb}^{(offline)}$ & $T_{rb}^{(online)}$ &
 $T_{fine}$ & $T_{fine} / T_{rb}$ & $L^2$ error  \\ [0.5ex] 
\hline\hline
$2 \times 10^{-4}$ & 52  & 6  & $157\second$ & $4.00\second$  & $4380\second$ & $27.2\times$ & 0.4997 \\
\hline
$1 \times 10^{-4}$ & 67  & 11 & $186\second$ & $6.09\second$  & $4380\second$ & $22.8\times$ & 0.1764 \\
\hline
$5 \times 10^{-5}$ & 77  & 15 & $206\second$ & $6.48\second$  & $4380\second$ & $20.6\times$ & 0.0648 \\
\hline
$2 \times 10^{-5}$ & 98  & 20 & $252\second$ & $7.51\second$  & $4380\second$ & $16.9\times$ & 0.0292 \\
\hline
$1 \times 10^{-5}$ & 115 & 26 & $283\second$ & $8.92\second$  & $4380\second$ & $15.0\times$ & 0.0136 \\
\hline
$5 \times 10^{-6}$ & 132 & 30 & $229\second$ & $10.03\second$ & $4380\second$ & $18.3\times$ & 0.0060 \\
\hline
$2 \times 10^{-6}$ & 150 & 38 & $254\second$ & $12.59\second$ & $4380\second$ & $16.4\times$ & 0.0038 \\
\hline
$1 \times 10^{-6}$ & 179 & 43 & $295\second$ & $15.77\second$ & $4380\second$ & $14.1\times$ & 0.0023 \\
\hline
$5 \times 10^{-7}$ & 194 & 48 & $315\second$ & $18.11\second$ & $4380\second$ & $13.2\times$ & 0.0018 \\
\hline
\end{tabular}                                                                            
\caption{\label{tab:result_bie}                                                               
Test results for our reduced basis method on the described Laplace equation problem. Here $\epsilon$
denotes the selection threshold used for reduced basis construction, $\skelCount$ denotes the number
of skeletons selected by our method (i.e., the number of fine solves used to construct the reduced
basis), $\reducedBasisDim$ denotes the dimension of the reduced basis constructed,
$T_{rb}^{(offline)}$ denotes the amount of time in seconds used in reduced basis construction,
$T_{rb}^{(online)}$ denotes the amount of time in seconds used to solve all problem instances from
$\paramSet$ using our method once the reduced basis has been constructed, $T_{rb} =
T_{rb}^{(offline)} + T_{rb}^{(online)}$ denotes the total amount of time in seconds used by our
method to compute approximations to all problem instances in $\paramSet$, and $T_{fine} / T_{rb}$
denotes the ratio between the time taken by our method to compute approximate solutions to all
problem instances in $\paramSet$ and the time $T_{fine}$ taken to naively compute all exact fine
solutions in $\paramSet$, i.e., the computational speed-up our algorithm provides. Finally, $L^2$
error denotes the average relative $L^2$ error, i.e., $\|\ve{u}(\param) -
\rbApproxParam{\param}\|_2 / \|\ve{u}(\param)\|_2$ averaged over the parameter set $\paramSet$.}
\end{table}

\begin{figure}[htb]
\centering
\includegraphics[width=0.5\textwidth]{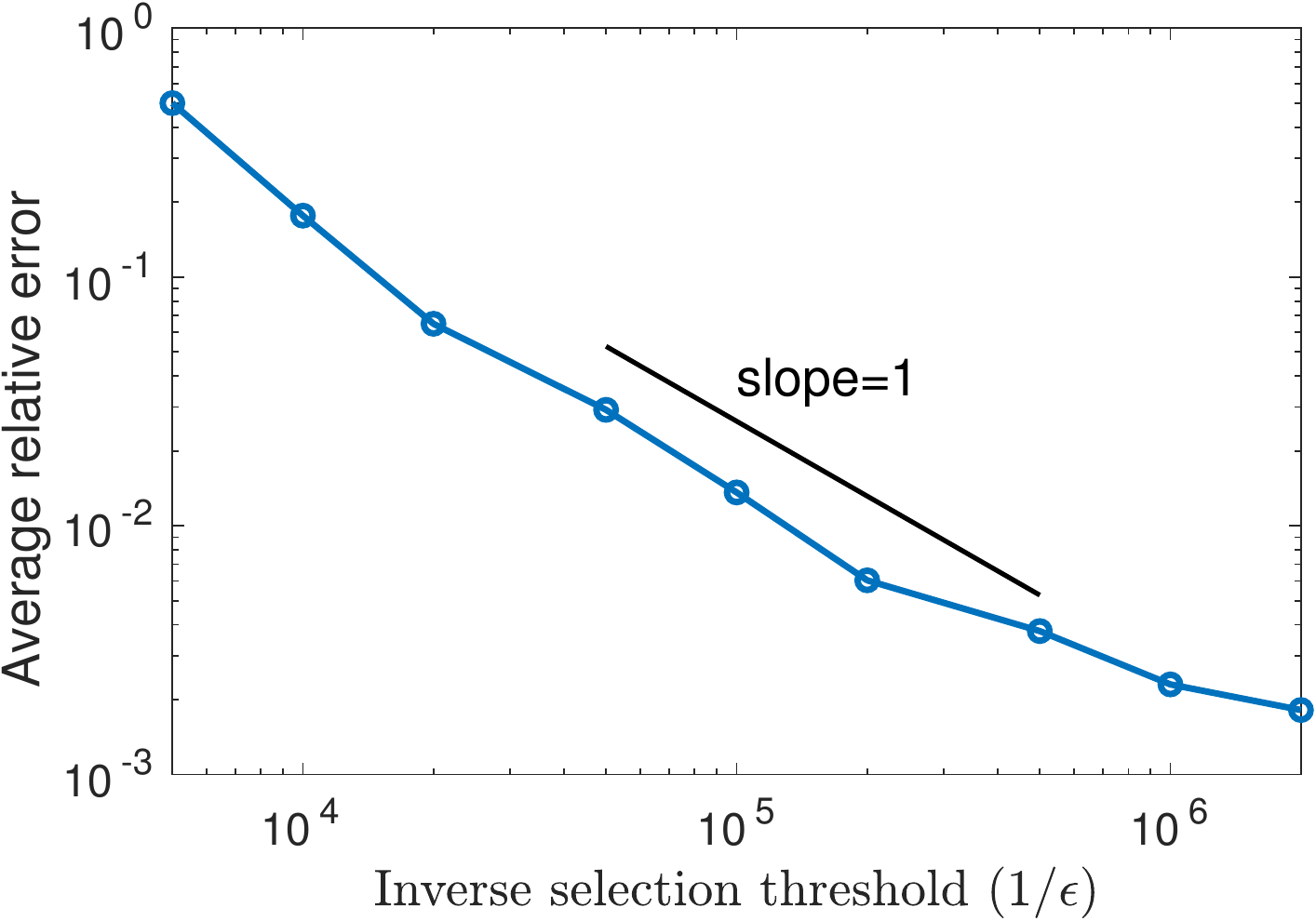}
\caption{\label{fig:conv_plot_bie} A log-log convergence plot of our method on the Laplace equation
example, showing the average error $\langle E \rangle$ over the parameter set $\paramSet$ against
the inverse of the selection threshold $\epsilon$.}
\end{figure}

Plot (a) in \Cref{fig:bie_Rsigma} presents the diagonal values of $R$ in \cref{basis_ext_alg} for this numerical experiment. As indicated in \cref{basis_ext_alg}, we compute these diagonal values from the coarse-proxy model solutions and use them to select skeleton parameters. Plot (b) in \Cref{fig:bie_Rsigma} shows the values of the
singular values $\sigma_k$ in \cref{mixing_matrix_alg}, which are used to construct the reduced basis.
One sees that, if the threshold $\epsilon$ satisfies $\epsilon<1\times 10^{-4}$, then $\log_{10}(R_{ii}/R_{11})$ and
$\log_{10}(\sigma_k/\sigma_1)$ decay almost linearly. We also provide a visualization of the first 36 reduced
basis (with respect to $\sigma_k$) in \cref{fig:basis_bie}. Note that, as the singular values $\sigma_k$ decay, the corresponding basis vectors contain increasingly more high
frequency information. For further evaluation, we provide \cref{fig:bie_examples}, a side-by-side comparison of our reduced basis method for $\epsilon=1\times 10^{-6}$ with the underlying full-order model on three different parameter instances. We note that our method provides an accurate approximation of the solution to the Laplace problem.

To quantitatively demonstrate the performance gains of our reduced-order model over the full-order model, as well as the reduced-order model's validity for different thresholds, we
evaluated both the full-order and reduced-order solutions for each element of $\paramSet$. The resulting relative $L^2$ error
$\|\ve{u}(\param) - \rbApproxParam{\param}\|_2 / \|\ve{u}(\param)\|_2$ and the runtime are all
recorded in \cref{tab:result_bie}. When evaluating all solutions in bulk, our reduced basis method provides between a ten-fold and twenty-fold performance increase over the full-order model for a wide range of accuracy targets between $.1\%$ and $10\%$ relative $L^2$ error, as recorded in \cref{bie_ie}. Note, crucially, that this figure \textit{includes the expensive offline phase of our reduced basis method}.
Moreover, if one halves the threshold $\epsilon$, the average relative $L^2$ shrinks by an approximate factor of two.  Therefore, this method approximately exhibits linear convergence with respect to the
threshold. This fact is further illustrated in \cref{fig:conv_plot_bie}.

We also note that the offline stage for $\epsilon = 5 \times 10^{-6}$ in \cref{tab:result_bie} takes
less time than the offline stage for $\epsilon = 1 \times 10^{-5}$, which is because proportionally
fewer skeletons are selected during the additional skeleton extraction phase.

\subsection{Radiative Transport Equation with Isotropic Scattering} \label{retsec}
We consider the steady state radiative transport equation with the form
\begin{equation} \label{rte_ie}
    \begin{split}
        \bm{v} \cdot \nabla_{\bm{x}} \Phi(\bm{x}, \bm{v}) + \mu_t(\bm{x}) \Phi(\bm{x}, \bm{v}) &= \frac{\mu_s(\bm{x})}{2\pi} \int_{\mathbb{S}^{1}} \Phi(\bm{x}, \bm{v}') \dd\bm{v}' + g(\bm{x}), \; \text{in } \domain \times \mathbb{S}^{1},\\
\Phi(\bm{x}, \bm{v}) &= 0, \; \text{on } \Gamma_{-},
    \end{split}
\end{equation}
where $\Phi(\bm{x}, \bm{v})$ denotes the photon flux at spatial position $\bm{x} \in \mathbb{R}^2$
in direction $\bm{v} \in \mathbb{S}^{1}$, and $g(\bm{x})$ is the light source. $\domain \subset \mathbb{R}^2$ is the problem domain, $\mathbb{S}^{1}$ is
the unit sphere in $\mathbb{R}^2$. $\Gamma_{-}$ is the inward facing problem boundary, given by
\begin{equation}
    \Gamma_{-} \equiv \{(\bm{x}, \bm{v}) \in \partial \domain \times \mathbb{S}^{1} \mid
    n(\bm{x}) \cdot v < 0 \},
\end{equation}
where $n(x)$ is the normal of domain $\domain$ at position $\bm{x}$. The boundary condition in
\eqref{retsec} enforces that no light is entering the domain of interest. 
The transport coefficient
$\mu_t(\bm{x})$ measures the total absorption at $\bm{x}$, which results from both physical
absorption as well as from scattering, the latter of which is quantified by the scattering coefficient
$\mu_s(\bm{x})$.

In this scenario, our quantity of interest is the local mean density
$m(\bm{x})$ defined as
\begin{equation}
  m(\bm{x}) \equiv \frac{1}{2\pi} \int_{\mathbb{S}^{1}} \Phi(\bm{x}, \bm{v}') \dd\bm{v}'.
\end{equation}
As studied in \cite{ren2016fast,fan2018fast}, one can reformulate the differential equation \cref{rte_ie} into an integral equation using the method of characteristics. This transformation
yields the integral equation
\begin{equation} \label{rte_ie2}
  \left[\frac{1}{\mu_s(\bm{x})} - \mathcal{K}\right]  u(\bm{x}) = \mathcal{K} g(\bm{x})
  \quad \text{ with }\quad
  \mathcal{K} \phi(\bm{x}) \equiv \int_{\domain} K(\bm{x}, \bm{y}) \phi(\bm{y}) \dd \bm{y},
\end{equation}
where $u(\bm{x}) = \mu_s(\bm{x}) m(\bm{x})$ 
the integral kernel $K(\bm{x}, \bm{y})$ of the operator $\mathcal{K}$ is given by
\begin{equation} \label{eq:int_kernel}
    K(\bm{x}, \bm{y}) \equiv \frac{1}{|\mathbb{S}^1|} \frac{1}{|\bm{x} - \bm{y}|} 
    \exp \left(- |\bm{x} - \bm{y}| \int_0^1 \mu_t(\bm{x} - \tau (\bm{x} - \bm{y}))\dd\tau\right).
\end{equation}

To bring this problem into the many-query setting of our reduced basis framework, we now suppose
that the scattering and transmission coefficients $\mu_s$ and $\mu_t$ have an explicit dependence on
a parameter $\param$ taken from some sample space $\paramSet_\infty$. Henceforth, we will therefore write them as
$\mu_s(\bm{x}; \param)$ and $\mu_t(\bm{x}; \param)$. Here, $\param$ can encode small fluctuations or
uncertainties about the underlying medium that the light propagates through. Making the dependence on the
parameter $\param$ explicit in \eqref{rte_ie2} gives us the set of integral equations to solve,
\begin{equation} \label{rte_final}
  \left[\mathcal{I} - \mu_s(\bm{x}; \param) \mathcal{K}(\param) \right] u(\bm{x}; \param) = \mu_s(\bm{x}; \param) \mathcal{K}(\param) g(\bm{x}) \,,
\end{equation}
where $\mathcal{I}$ is the identity operator. 

To discretize the above equation, we use a collocation method combined with Gauss-Legendre quadrature, as outlined in \cite{fan2018fast}. This discretization gives us the linear system
\begin{equation} \label{solvs}
    \operator{\param} \solution{\param} = \mathsf{f}(\param) \,,
\end{equation}
where $\operator{\param}$ and $\mathsf{f}(\param)$ have the forms
\begin{equation}
\operatorNP(\param) = \mat{I} + \operatorAffineParam{\param}, \qquad \source \equiv -\operatorAffineParam{\param} \ve{g} \,,
\end{equation}
and $\mat{I}$, $\operatorAffineParam{\param}$, and $\ve{g}$ are the discretized versions of $\mathcal{I}$, $- \mu_s(\bm{x}; \param) \mathcal{K}(\param)$, and $g(\bm{x})$ respectively.

The application of our framework to this problem is now straightforward. To solve the full-order
model, we use hierarchical interpolative factorization \cite{l2016hierarchical}. For our coarse-proxy
model, we simply use significantly fewer collocation points in $\domain$. We then use the method
described in \cref{basis_ext_sec,rb_construct} to compute a suitable reduced basis matrix $\reducedBasis$ for this problem in the offline stage.

In the remainder of the offline stage, we use the method described in \cref{offset_ops,redops} to sample $\operatorAffineParam{\param}$ and construct a mixing matrix $\mathsf{M}$ such that
\begin{equation} \label{eq:rb_op_construct1}
    \mathcal{B}_{rb} \approx \skel{\mathcal{B}}_{rb} \mathsf{M} ,
\end{equation}
where once again we define $\mathcal{B}_{rb}$ to be the matrix vectorized reduced operators,
\begin{equation}
  \mathcal{B}_{rb} \equiv \begin{bmatrix}
    \text{vec}(\reducedBasis^\T \operatorAffineParam{\param_1} \reducedBasis) &
    \text{vec}(\reducedBasis^\T \operatorAffineParam{\param_2} \reducedBasis) & \dots &
    \text{vec}(\reducedBasis^\T \operatorAffineParam{\param_n} \reducedBasis)
  \end{bmatrix},
\end{equation}
and $\skel{\mathcal{B}}_{rb} \equiv \mathcal{B}_{rb}(:, \mathfrak{S})$ are our reduced operator skeletons. For the samples $\mathfrak{O}$ in the computation of $\mathsf{M}$, we use a small number of randomly selected columns.

With the offline stage complete, we now switch focus to the online stage. To solve the desired equation
\begin{equation}
  (\mathsf{I} + \reducedBasis^\T \operatorAffineParam{\param} \reducedBasis) \rbCoeffParam{\param}
  = - \reducedBasis^\T \operatorAffineParam{\param} \mathsf{g} ,
  \quad \solution{\param} \approx \reducedBasis \, \rbCoeffParam{\param},
\end{equation}
we can assemble the reduced operator $\mathsf{I} + \reducedBasis^\T \operatorAffineParam{\param}
\reducedBasis$ from the reduced operator skeletons $\mathcal{B}_{rb}$ by using \eqref{eq:rb_op_construct1}. Consider the matrix
\begin{equation}
  \mathcal{F}_{rb} \equiv 
  \begin{bmatrix}
    -\reducedBasis^\T \operatorAffineParam{\param_1} \mathsf{g} & 
    -\reducedBasis^\T \operatorAffineParam{\param_2} \mathsf{g} & \dots & 
    -\reducedBasis^\T \operatorAffineParam{\param_n} \mathsf{g} 
  \end{bmatrix},
\end{equation}
and note that the interpolation weights computed for the reduced operators $\reducedBasis^\T \operatorAffineParam{\param} \reducedBasis$ carry over to this matrix. That is,
\begin{equation}
  \mathcal{F}_{rb} \approx \skel{\mathcal{F}}_{rb} \mathsf{M} ,
\end{equation}
where $\skel{\mathcal{F}}_{rb} \equiv \mathcal{F}(:, \mathfrak{S})$. Since we must assemble the quantities
$-\operatorAffineParam{\skel{\param}_i} \mathsf{g}$ during the computation of
fine solutions for our reduced basis regardless, computing the matrix
$\skel{\mathcal{F}}_{rb}$ is fairly inexpensive, and only involves applying the matrix
$\reducedBasis^\T$ to the vectors $-\operatorAffineParam{\skel{\param}_i} \mathsf{g}$. This means, if we compute $\skel{\mathcal{F}}_{rb}$ during the offline stage, we then have an inexpensive method of assembling the right hand side, regardless of the fact that the expression involves the operator $\operatorAffineParam{\param}$.

\subsubsection{Results}
As a test case for the above example, we consider the domain $\domain \equiv [0, 1]^2$. We let $\mu_s(\bm{x}; \param)$ and $\mu_t(\bm{x}; \param)$ be Guassians with varying centers and widths,
\begin{align}
\mu_t(\bm{x}; \param) &\equiv \mu_s(\bm{x}; \param) \equiv 1 + A_\param \exp(-((x_1 - c_{1\param})^2 - (x_2 - c_{2\param})^2) / \theta_\param^2) ,
\end{align}
where the parameters $\param \in \paramSet$ have the form
\begin{equation}
  \param \equiv \begin{bmatrix} A_\param & c_{1\param} & c_{2\param} & \theta_\omega \end{bmatrix}.
\end{equation}
We take the source term $g(\bm{x})$ to be
\begin{equation}
g(\bm{x}) \equiv \exp(-256((x_1-0.5)^2+(x_2-0.5)^2)) .
\end{equation}
To build the parameter space $\paramSet$ we vary both the width and the location of the Gaussian ensemble above. Let $\paramSet_{A, \theta, N}$ be defined as
\begin{equation}
  \paramSet_{A, \theta, N} \equiv \left\{ \begin{bmatrix} A & i / N & j / N& \theta \end{bmatrix} \mid i, j = 0, ..., N \right\},
\end{equation}
that is, parameters for Gaussians with width $\theta$ and amplitude $A$ centered at grid points $(i/N, j/N)$. Take our parameter space $\paramSet$ to consist of these Gaussians with three different widths/amplitudes,
\begin{equation}
  \paramSet \equiv \bigcup_{A \in \mathcal{A}} \bigcup_{\theta \in \Theta} \paramSet_{A, \theta, N},
\end{equation}
where
\begin{equation*}
\mathcal{A} \equiv \{ 2, 4, 6, 8, 10 \}, \qquad \Theta \equiv \{ 0.2, 0.3, 0.4, 0.5, 0.6 \},
\quad N=20.
\end{equation*}
This gives a total parameter space size of $|\paramSet| = 11025$. Our full-order model is the model
described in \cref{retsec}, with a grid size of $n_f \times n_f$ where $n_f = 128$. We use the
algorithm described in this paper to build a reduced basis for this model and approximate true
solutions. For our coarse-proxy model with the same model described in \cref{retsec} but with a grid
size of $n_c \times n_c$ instead, where $n_c = 32$. We use the procedure described in \cref{add_sk}
to add additional skeletons to our skeleton set when operator samples cannot be well-represented
using a linear combination the operator samples of the selected skeletons. For selecting additional
skeletons, we use a selection threshold multiplier of $\epsilonMul = 1.5$. 

\begin{figure}[htb]
 \centering 
 \subfloat[Diagonals of $R$]{
   \includegraphics[width=0.4\textwidth]{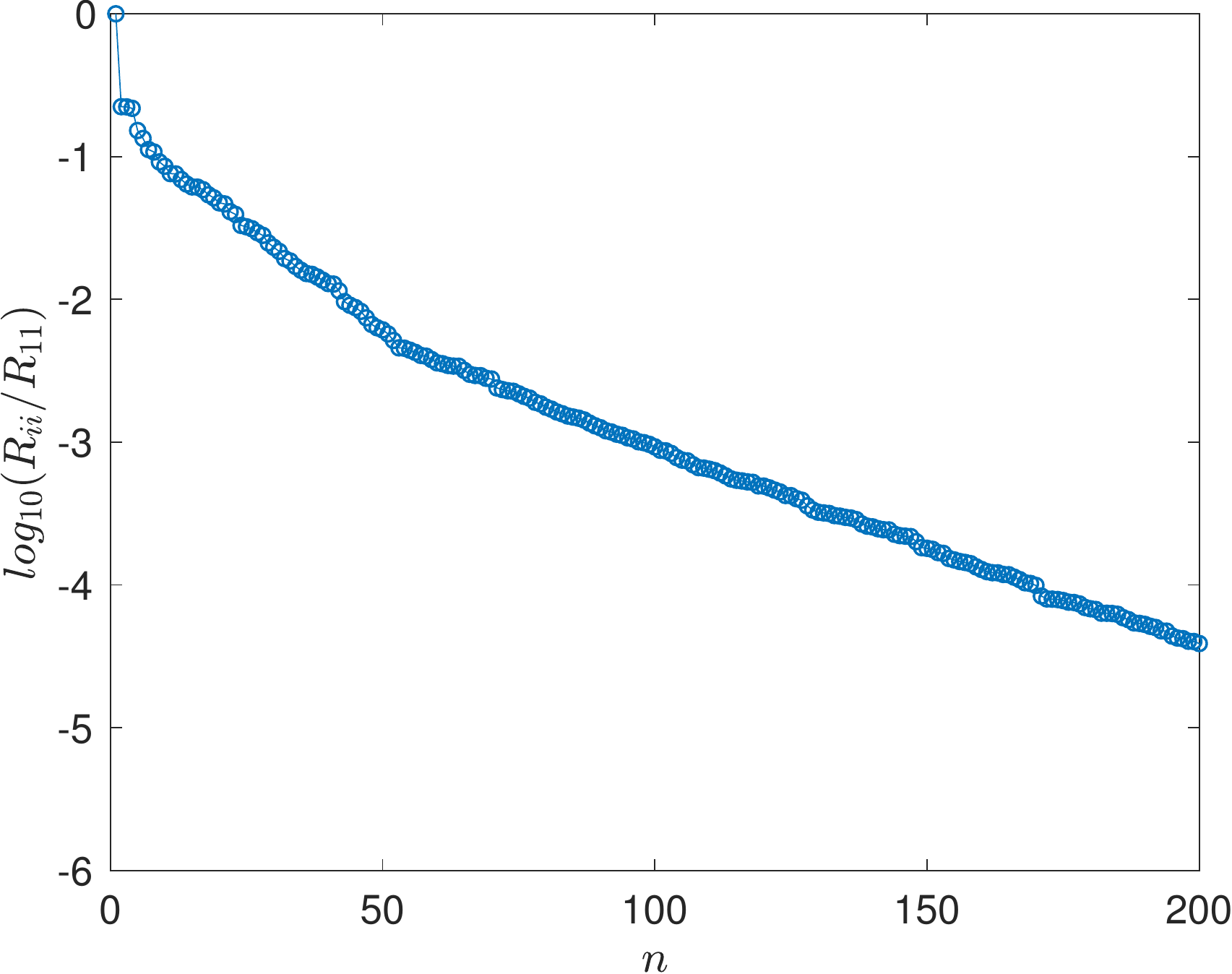}
 }\hspace{0.1\textwidth}
 \subfloat[Sigular values $\sigma_k$]{
   \includegraphics[width=0.4\textwidth]{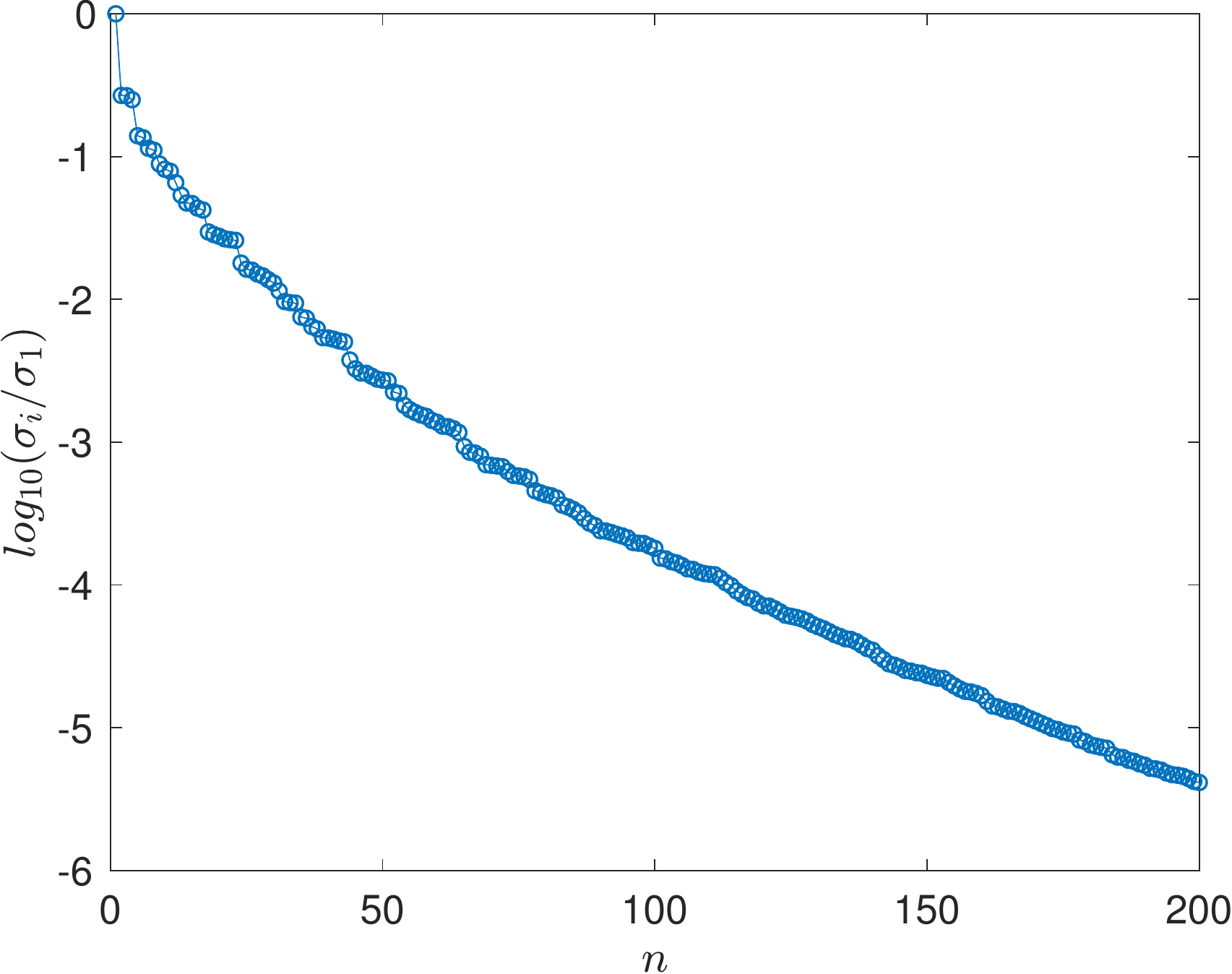}
 }
 \caption{\label{fig:rte_Rsigma} Left: the normalized diagonal values of $R$ in the skeleton selection for the described radiative transport problem. Right:
 the normalized singular values in the basis construction for the described radiative transport problem.}
\end{figure}

\begin{figure}[htb]
  \centering
  \includegraphics[width=\textwidth]{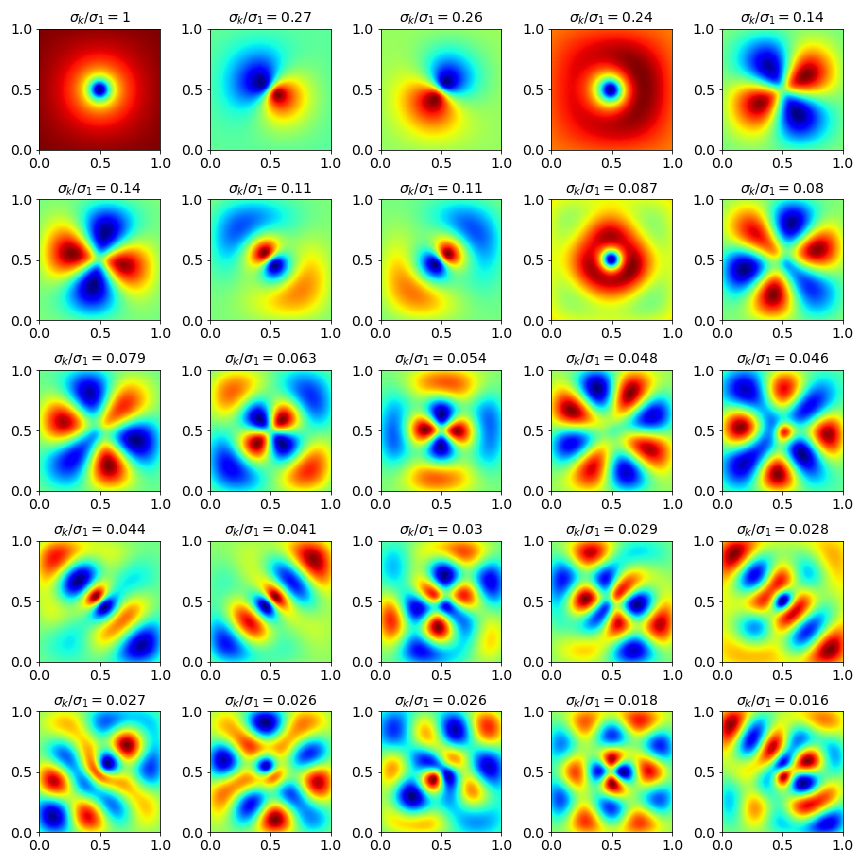}
      \caption{ \label{fig:basis_rte} Reduced basis vectors generated for the described
    radiative transport problem and their corresponding normalized singular values $\tilde{\sigma}_k$.}
\end{figure}

Plot (a) in \Cref{fig:rte_Rsigma} presents the diagonal values of $R$ in \cref{basis_ext_alg} for this numerical experiment. As indicated in \cref{basis_ext_alg}, we compute these diagonal values from the coarse-proxy model solutions and use them to select skeleton parameters. Plot (b) in \Cref{fig:rte_Rsigma} shows the values of the
singular values $\sigma_k$ in \cref{mixing_matrix_alg}, which are used to construct the reduced basis for $\epsilon = 5 \times 10^{-6}$. We also provide a visualization of the first 25 reduced
basis (with respect to $\sigma_k$) in \cref{fig:basis_rte}. Note that, as the singular values $\sigma_k$ decay, the corresponding basis vectors contain increasingly more high
frequency information.

\begin{figure}[htb]
  \centering
  \includegraphics[width=0.99\textwidth]{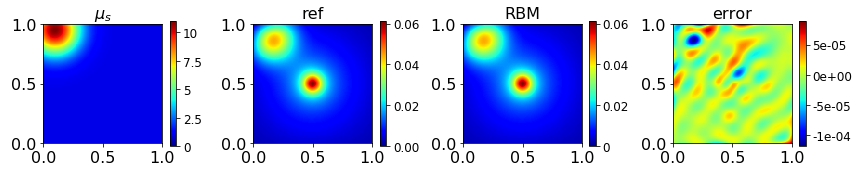}
  \includegraphics[width=0.99\textwidth]{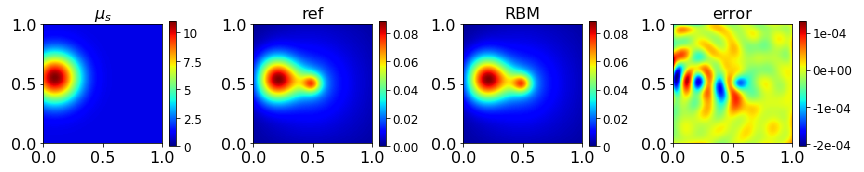}
  \includegraphics[width=0.99\textwidth]{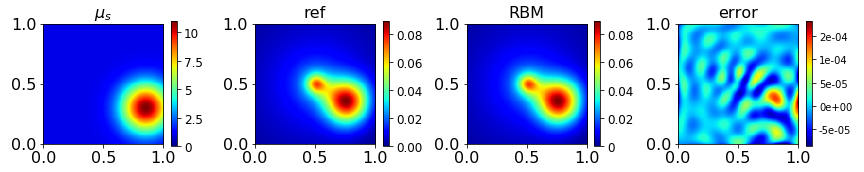}
  \caption{\label{fig:rte_examples} Three examples for the solutions evaluated by the reduced basis
    method for the radiative transport equation for the parameter set $\Omega$ with threshold
  $\epsilon=1\times 10^{-4}$ and their corresponding reference solutions and error. }
\end{figure}

\begin{table}[htb]
\centering
 \begin{tabular}{|| c || c | c || c | c| c|c || c ||} 
 \hline
 $\epsilon$ & $\skelCount$ & $\reducedBasisDim$ & $T_{rb}^{(offline)}$ & $T_{rb}^{(online)}$ &
 $T_{fine}$ & $T_{fine} / T_{rb}$ & $L^2$ error \\ [0.5ex] 
 \hline\hline                                                                                               
$ 1 \times 10^{-3}$ & 97 & 63 & $1274 \, \text{sec}$ & $12.25 \, \text{sec}$ & $46676 \, \text{sec}$ & $36.3 \times$ & 0.0552 \\                                     
\hline                                                                                               
$ 5 \times 10^{-4}$ & 118 & 75 & $1478 \, \text{sec}$ & $21.97 \, \text{sec}$ & $46676 \, \text{sec}$ & $31.1 \times $ & 0.0491 \\                                    
\hline                                                                                               
$2 \times 10^{-4}$ & 148 & 95 & $1489 \, \text{sec}$ & $83.64 \, \text{sec}$ & $46676 \, \text{sec}$ & $29.7 \times$ & 0.0406 \\                                    
\hline                                                                                               
$1 \times 10^{-4}$ & 169 & 111 & $1565 \, \text{sec}$ & $110.2 \, \text{sec}$ & $46676 \, \text{sec}$ & $27.9\times$ & 0.0307 \\                                  
\hline                                                                                               
$5 \times 10^{-5}$ & 192 & 127 & $1786 \, \text{sec}$ & $144.7\, \text{sec}$ & $46676 \, \text{sec}$ & $24.2 \times $ & 0.0249 \\                                  
\hline                                                                                               
$2 \times 10^{-5}$ & 225 & 153 & $1995 \, \text{sec}$ & $187.51$ & $46676 \, \text{sec}$& $21.4\times$ & 0.0135 \\                                  
\hline                                                                                               
$1 \times 10^{-5}$ & 250 & 170 & $2095 \, \text{sec}$ & $239.1 \, \text{sec}$ & $46676 \, \text{sec}$ & $20.0\times$ & 0.0096 \\                                  
\hline                                                                                               
$5 \times 10^{-6}$ & 277 & 193 & $2388 \, \text{sec}$ & $303.8 \, \text{sec}$ & $46676 \, \text{sec}$ & $17.3 \times$ & 0.0081 \\                                  
\hline
\end{tabular} \\[10pt]
\caption{ \label{tab:result_rte}
  Test results for the reduced basis method on the described radiative transport problem. Here
  $\epsilon$ denotes the selection threshold used for reduced basis construction, $\skelCount$ denotes
  the number of skeletons selected by our method (i.e., the number of fine solves used to construct the
  reduced basis), $\reducedBasisDim$ denotes the dimension of the reduced basis constructed,
  $T_{rb}^{(offline)}$ denotes the amount of time in seconds used in reduced basis construction,
  $T_{rb}^{(online)}$ denotes the amount of time in seconds used to solve all problem instances from
  $\paramSet$ using our method once the reduced basis has been constructed, $T_{rb} =
  T_{rb}^{(offline)} + T_{rb}^{(online)}$ denotes the total amount of time in seconds used by our
  method to compute approximations to all problem instances in $\paramSet$, and $T_{fine} / T_{rb}$
  denotes the ratio between the time taken by our method to compute approximate solutions to all
  problem instances in $\paramSet$ and the time $T_{fine}$ taken to naively compute all exact fine
  solutions in $\paramSet$, i.e., the computational speed-up our algorithm provides. Finally, 
  $L^2$ error denotes the average relative $L^2$ error, i.e., $\|\ve{u}(\param) -
  \rbApproxParam{\param}\|_2 / \|\ve{u}(\param)\|_2$ averaged over the parameter set $\paramSet$. 
}
\end{table}

Once again for the visualization proposes, we present the reduced basis for $\epsilon=1 \times 10^{-4}$ in \cref{fig:basis_rte} --- together with side-by-side comparisons, on three different parameter instances, of our reduced basis approximation for $\epsilon=1 \times 10^{-4}$ to the underlying full-order solution, in \cref{fig:rte_examples}. We note the reduced basis method gives a good
approximation of the solution.

To further test the validity and efficiency of our reduced-order model, we run a parallel battery of tests to those we ran for the previous numerical example.
For each element of $\paramSet$, we compute both the true solution $\ve{u}(\param)$ and the reduced basis approximation
$\rbApproxParam{\param}$, and afterwards evaluate the relative $L^2$ error $\|\ve{u}(\param) - \rbApproxParam{\param}\|_2
/ \|\ve{u}(\param)\|_2$. 
We present the results of this computation in \cref{tab:result_rte}. 
When evaluating all solutions in bulk, our reduced basis method provides between a seventeen-fold and thirty-five-fold performance increase over the full-order model for a wide range of accuracy targets between $.8\%$ and $5\%$ relative $L^2$ error, as recorded in \cref{bie_ie}. Note, once again, that this figure \textit{includes the expensive offline phase of our reduced basis method}.
While this numerical example does not quite match the linear convergence of the previous numerical example, we still note that the error always decreases as the parameter $\epsilon$ decreases. Therefore, our method exhibits convergence, as seen in \cref{tab:result_rte}.

\section{Conclusion and Future Work}
We have developed a simple and general-purpose reduced basis approximation technique for linear elliptic
integral operators. As shown by the empirical results, this method results in significant
performance increases on the simple problems we have applied it to.  Due to the complexity scaling
exhibited by numerical simulations, this method might produce even more significant performance
increases at scale.  Moreover, we hope that the techniques put forth in this paper will provide a
useful starting point for future work in model order reduction for integral equations.

Possible avenues for such future work include the application of these techniques to larger scale
problems, or perhaps, the application of these techniques to electromagnetic scattering to give a
real comparison to currently existing work in \cite{fares2011reduced, hesthaven2012certified}. Other
possible areas for future work include the method by which the reduced operators $\reducedBasis^\T
\operator{\param} \reducedBasis$ are assembled. One could imagine finding a better operator sampling
mask than the randomly selected ones in this paper. It may also be possible that at scale, the
method we use to compute interpolation coefficients may be subject to overfitting. However, in our
experience working on the radiative transport and Laplace equation examples, the interpolation error
for the reduced operators is not a dominant source of error except at very small values of the
threshold $\epsilon$ where the total average error is already very small. Finally, for problems at
scale, there is a trade-off that must be made between the quality and computation time for the
coarse-proxy model. It may be useful in this situation to use a series of coarse-proxy models (rather than a
single one), each subsequent one finer than the last, to progressively filter down the parameter set
$\paramSet$ to the skeleton set $\skel{\paramSet}$. 

\bibliographystyle{abbrv}
\bibliography{references}

\end{document}